\documentclass[11pt]{amsart}
\usepackage{stmaryrd}
\usepackage{amsmath,amssymb,amsthm}
\usepackage[pagebackref]{hyperref}
\hypersetup{colorlinks=true,urlcolor=blue,citecolor=blue,linkcolor=blue}

\usepackage{tikz}
\usepackage{tikz-cd}
\usepackage{color}
\usepackage{enumerate}
\usepackage{booktabs}
\usepackage{colonequals}
\usepackage{longtable}
\usepackage{caption}
\usepackage{float}

\usepackage[protrusion=false]{microtype}

\usepackage{arydshln}

\usepackage[symbol]{footmisc}

\def\ol#1{\overline{#1}}% 		overline
\theoremstyle{plain}
    \newtheorem{theorem}{Theorem}[section]

    \newtheorem{proposition-definition}[theorem]{Proposition-Definition}
    \newtheorem{lemma-definition}[theorem]{Lemma-Definition}
\theoremstyle{definition}

    \newtheorem{remark}[theorem]{Remark}

\newtheoremstyle{named}{}{}{\itshape}{}{\bfseries}{.}{.5em}{\thmnote{#3 }#1}
\theoremstyle{named}
\newtheorem*{namedtheorem}{Theorem}
   
\usepackage{cleveref}
\crefname{theorem}{Theorem}{Theorems}
\crefname{lemma}{Lemma}{Lemmata}
\crefname{corollary}{Corollary}{Corollaries}
\crefname{proposition}{Proposition}{Propositions}
\crefname{proposition-definition}{Proposition-Definition}{Proposition-Definitions}
\crefname{lemma-definition}{Lemma-Definition}{Lemma-Definitions}
\crefname{definition}{Definition}{Definitions}
\crefname{conjecture}{Conjecture}{Conjectures}
\crefname{question}{Question}{Questions}
\crefname{example}{Example}{Examples}
\crefname{algorithm}{Algorithm}{Algorithms}
\crefname{remark}{Remark}{Remarks}
\crefname{assumption}{Assumption}{Assumptions}
    
\def\Alphabet{A,B,C,D,E,F,G,H,I,J,K,L,M,N,O,P,Q,R,S,T,U,V,W,X,Y,Z}%  Capitalized Alphabet
\def\alphabet{a,b,c,d,e,f,g,h,i,j,k,l,m,n,o,p,q,r,s,t,u,v,w,x,y,z}%	lowercase alphabet
\def\endpiece{xxx}%									marks end of list
\def\makeAlphabet[#1]{\expandafter\makeA#1,xxx,}%		Ex. \makeAlphabet[A,B]
\def\makealphabet[#1]{\expandafter\makea#1,xxx,}%		Ex. \makealphabet[c,d]
\def\makeA#1,{\def\temp{#1}\ifx\temp\endpiece\else%
\mkbb{#1}\mkfrak{#1}\mkbf{#1}\mkcal{#1}\mkscr{#1}\mkbs{#1}\expandafter\makeA\fi}%
\def\makea#1,{\def\temp{#1}\ifx\temp\endpiece\else\mkfrak{#1}\mkbf{#1}\mkbs{#1}\expandafter\makea\fi}%
\def\mkbb#1{\expandafter\def\csname bb#1\endcsname{\mathbb{#1}}}%      Define bb
\def\mkfrak#1{\expandafter\def\csname fr#1\endcsname{\mathfrak{#1}}}%    Define frak
\def\mkbf#1{\expandafter\def\csname b#1\endcsname{\mathbf{#1}}}%           Define bold letters
\def\mkcal#1{\expandafter\def\csname c#1\endcsname{\mathcal{#1}}}%       Define calligraphy
\def\mkscr#1{\expandafter\def\csname s#1\endcsname{\mathscr{#1}}}%       Define script
\def\mkbs#1{\expandafter\def\csname bs#1\endcsname{{\boldsymbol{#1}}}}%       Define bold symbol
\def\makeop[#1]{\xmakeop#1,xxx,}%					Ex. \makeop[Hom,Spec]
\def\mkop#1{\expandafter\def\csname #1\endcsname{{\mathrm{#1}}}} % 
\def\xmakeop#1,{\def\temp{#1}\ifx\temp\endpiece\else\mkop{#1}\expandafter\xmakeop\fi}%
\def\makeup[#1]{\xmakeup#1,xxx,}%					Ex. \makeop[Hom,Spec]
\def\mkup#1{\expandafter\def\csname #1\endcsname{{\mathrm{#1}\,}}} % 
\def\xmakeup#1,{\def\temp{#1}\ifx\temp\endpiece\else\mkup{#1}\expandafter\xmakeup\fi}%
\makeAlphabet[\Alphabet]%				Define bb, frak, bf, cal for Capitalized Alphabet
\makealphabet[\alphabet]%  				Define frak and bf for uncapitalized alphabet
\makeop[Hom,Tor,Sel,GL,H,R,ord,Pic,NS,Br,BR,Zar,fppf,rig,cris,syn,dR,Gal,Ker,Im,Coker,ker,coker,BSD,Res,T,Frob,rk,sh,nr,Re,Im,Jac,Aut,End,gcd]
\makeup[Spec,Proj,Spwf,Sp,Sh,Spf,Sch,id,Tot,sm,an,tor,holim,hocolim,dim,red,rank,rk,CH,Tr,GM]
  \makeatletter
    
    \@addtoreset{equation}{section}
  \makeatother

\makeatletter
\newcommand{\labitem}[2]{%
\def\@itemlabel{\textbf{#1}}
\item
\def\@currentlabel{#1}\label{#2}}
\makeatother

\newcommand*{\sheafhom}{\sH\kern-.5pt \mathrm{om}}
\newcommand*{\sheafext}{\sE\kern-.5pt \mathrm{xt}}

\newcommand{\Q}{\bbQ}
\newcommand{\Z}{\bbZ}
\newcommand{\F}{\bbF}

\renewcommand{\phi}{\varphi}
\renewcommand{\epsilon}{\varepsilon}
\renewcommand{\theta}{\vartheta}

\newcommand{\iso}{\simeq}

\newcommand{\isoto}{\xrightarrow{\sim}}

\newcommand{\inj}{\hookrightarrow}

\DeclareFontFamily{U}{wncy}{}
\DeclareFontShape{U}{wncy}{m}{n}{<->wncyr10}{}
\DeclareSymbolFont{mcy}{U}{wncy}{m}{n}
\DeclareMathSymbol{\Sha}{\mathord}{mcy}{"58}

\definecolor{darkred}{rgb}{0.8,0,0}

\usepackage{comment}

\begin{document}

\title[Rational points on hyperelliptic Atkin-Lehner quotients]{Rational points on hyperelliptic\\ Atkin-Lehner quotients of modular curves\\ and their coverings}

\author[N. Adžaga]{Nikola Adžaga}
\address{Department of Mathematics, Faculty of Civil Engineering, University of Zagreb}
\email{nadzaga@grad.hr}

\author[S. Chidambaram]{Shiva Chidambaram}
\address{Department of Mathematics, MIT}
\email{shivac@mit.edu}

\author[T. Keller]{Timo Keller}
\address{Lehrstuhl Mathematik~II (Computeralgebra)\\Universität Bayreuth\\Universitätsstraße 30\\95440 Bayreuth, Germany} 
\curraddr{Leibniz Universität Hannover, Institut für Algebra, Zahlentheorie und Diskrete Mathematik, Welfengarten 1, 30167 Hannover, Germany}
\email{keller@math.uni-hannover.de}
\urladdr{\url{https://www.timo-keller.de}}

\author[O. Padurariu]{Oana Padurariu}
\address{Department of Mathematics and Statistics, Boston University}
\email{oana@bu.edu}

\thanks{N.\@ A.\@ is supported by the Croatian Science Foundation under the project no. IP2018-01-1313.
S.\@ C.\@ is supported by the Simons Foundation grant \#550033.
T.\@ K.\@ is supported by the Deutsche Forschungsgemeinschaft (DFG), Projektnummer STO 299/18-1, AOBJ: 667349 while working on this article.
O.\@ P.\@ is supported by NSF grant DMS-1945452 and Simons Foundation grant \#550023.
}

\date{\today}
\maketitle

\begin{abstract}
    We complete the computation of all $\Q$-rational points on all the $64$ maximal Atkin-Lehner quotients $X_0(N)^*$ such that the quotient is hyperelliptic. To achieve this, we use a combination of various methods, namely the classical Chabauty--Coleman, elliptic curve Chabauty, quadratic Chabauty, and the bielliptic quadratic Chabauty method (from a forthcoming preprint of the fourth-named author) combined with the Mordell-Weil sieve. Additionally, for square-free levels $N$, we classify all $\Q$-rational points as cusps, CM points (including their CM field and $j$-invariants) and exceptional ones. We further indicate how to use this to compute the $\Q$-rational points on all of their modular coverings.
\end{abstract}

\section{Introduction}
Let $N$ be a positive integer and let $X_0(N)$ be the smooth compactification of the affine modular curve $Y_0(N)$ of level $N$. It is the moduli space parametrizing cyclic $N$-isogenies between (generalized) elliptic curves. For every $d$ dividing $N$ such that $\gcd(d, N/d) = 1$, let $w_d$ denote the corresponding Atkin-Lehner involution on $X_0(N)$. As explained in \cite[\S{I.2}]{GrossZagier1986}, the involution $w_d$ sends a point $x \in Y_0(N)$ representing the isogeny $\varphi \colon E \to E'$, to the point representing the composite isogeny
\[
    E/\ker(\phi)[d] \to E/\ker(\phi) \isoto E' \to E'/\ker(\phi^\vee)[d].
\]
By the valuative criterion for properness, the involution $w_d$ has a unique extension to $X_0(N)$.

Let $W(N)$ denote the group generated by all the Atkin-Lehner involutions on $X_0(N)$. Then $W(N)$ is isomorphic to $(\Z/2\Z)^{\omega(N)}$ where $\omega(N)$ is the number of distinct prime factors of $N$. The main object of study in this paper is the quotient curve $X_0(N)^* \colonequals X_0(N)/W(N)$.
We denote its Jacobian by $J_0(N)^*$. We also use the notation $X_0(N)^+ \colonequals X_0(N)/w_N$. (Note that $w_N$ sends an isogeny to its dual.)

If a rational point on $X_0(N)^*$ is the image of a cusp on $X_0(N)$, we call it a cuspidal point or simply a cusp. For every non-cuspidal rational point on $X_0(N)^*$, its lifts in $X_0(N)$ correspond to $\Q$-curves defined over multi-quadratic extensions of $\Q$~\cite{Elkies2004}. Recall that a $\Q$-curve is an elliptic curve defined over $\overline{\Q}$ that is isogenous to all of its Galois conjugates. If the $\Q$-curves have complex multiplication, we will call the rational point a CM point. Otherwise, we call it an exceptional point. Since every CM elliptic curve is a $\Q$-curve, the most interesting rational points on $X_0(N)^*$ are the exceptional points corresponding to non-CM $\Q$-curves.

Mazur~\cite{MazurEisenstein} called the problem of determining the rational points of Atkin-Lehner quotients of $X_0(N)$ an ``extremely interesting diophantine question''. For example, $\Q$-rational points on $X_0(N)^+$ lift to points of degree $1$ or $2$ on $X_0(N)$; for the computation of quadratic points on $X_0(N)$, see~\cite{Box2021}.

When $X_0(N)^*$ has genus $0$ or $1$, the sets $X_0(N)^*(\Q)$ of rational points were determined in~\cite{GL1998} who also gave explicit parametrization of the families of $\Q$-curves described by these points. When $N$ is a prime power, $X_0(N)^* = X_0(N)^+$, since $w_N$ is the only non-trivial Atkin-Lehner involution. For $N$ prime and $\text{genus}(X_0(N)^*) \in \{2,3\}$, the sets $X_0(N)^*(\Q)$ were computed in~\cite{BDMTV21, BBBLMTV21}. For $N$ prime and $\text{genus}(X_0(N)^*) \in \{4, 5, 6\}$, the same was carried out by the first three authors and Vishal Arul, Lea Beneish, Mingjie Chen, and Boya Wen in~\cite{AABCCKW}. For $N$ composite and $\text{genus}(X_0(N)^+) \in \{2, 3, 4, 5, 6\}$, the articles~\cite{Momose87, Galbraith2002, AraiMomose10, ArulMueller} complete the determination of the sets $X_0(N)^+(\Q)$. These articles together verify a conjecture of Galbraith~\cite{Galbraith2002} on the set of exceptional rational points on the curves $X_0(N)^+$ with genus between $2$ and $5$. Together with these articles and the work of~\cite{BGX21}, in this paper, we finish the determination of rational points on all hyperelliptic curves $X_0(N)^*$.

The class of $\Q$-curves is important because it is a mild generalization of the class of elliptic curves over $\Q$~\cite{Ell2004}. Analogous to the compatible system of $\ell$-adic Galois representations associated to the Tate module of an elliptic curve, one naturally associates to a $\Q$-curve $E$, a compatible system of $\ell$-adic Galois representations
\[ \rho_{E,\ell} : \Gal(\overline{\Q}/\Q) \rightarrow \overline{\Q}^{\times}\GL_2(\Q_{\ell}). \]
Ribet~\cite{Rib2004} conjectured that $\Q$-curves are exactly all the elliptic curves over $\overline{\Q}$ which are quotients of the modular curves $X_1(N)$. In other words, if $E$ is a $\Q$-curve, then the associated representation $\rho_{E,\ell}$ is the $\ell$-adic Galois representation attached to a cusp form on $\Gamma_1(N)$ for some $N$. Several partial results~\cite{HHM99, Hid2000, ES01} for many classes of $\Q$-curves were obtained using Wiles' modularity lifting theorems, until the proof of Serre's conjecture by Khare, Kisin, and Wintenberger~\cite{KW2009} proved Ribet's conjecture completely~\cite[Theorem~7.2]{Kha2010}.

Finding the set $X_0(N)^*(\Q)$ is an arithmetically interesting problem for any level $N$. But for the purpose of studying non-CM $\Q$-curves up to isogeny, the following result of Elkies says that it is enough to look at exceptional rational points on $X_0(N)^*$ for squarefree levels $N$.

\begin{theorem}[{\cite{Elkies2004}}]
Let $E$ be a $\Q$-curve without complex multiplication. Then there exists some square-free $N > 0$ such that $E$ is geometrically isogenous to the $\Q$-curves associated to an exceptional rational point on $X_0(N)^*$.
\end{theorem}

The following result due to Hasegawa classifies all the hyperelliptic Atkin-Lehner quotients $X_0(N)^*$.

\begin{theorem}[{\cite[Theorem~B]{Hasegawa97}}] \label[theorem]{hyperelliptic star quotients}
There are $64$ values of $N$ for which $X_0(N)^*$ is hyperelliptic. Of these, there are only $7$ values of $N$ for which $X_0(N)^*$ is hyperelliptic with genus $g \geq 3$, namely $N = 136$, $171$, $207$, $252$, $315$ ($g = 3$), $176$ ($g = 4$), and $279$ ($g = 5$).
\end{theorem}

Among the $57$ levels $N$ for which $X_0(N)^*$ has genus $2$, there are $39$ square-free levels\footnote[2]{There seems to be a typo in~\cite{BGX21} regarding their number. They say that there are $36$ square-free levels. They miss the levels $166, 255, 330$.} and $18$ non-squarefree levels. The authors in~\cite{BGX21} compute for a subset of these square-free levels the set of rational points on $X_0(N)^*$ using elliptic curve Chabauty. They also compute the $\Q$-curves corresponding to each non-cuspidal point in the set. The prime levels $N = 67$, $73$, $103$, $107$, $161$, $167$, $191$ were approached using quadratic Chabauty in~\cite{BBBLMTV21} and~\cite{BDMTV21}. We compute $X_0(N)^*(\Q)$ for the remaining square-free levels
\begin{equation}
\label{geq2_sqfree_levels}
    \{ 133,134,146,166\footnote[4]{citing forthcoming work~\cite{BianchiPadurariu}; note that our forthcoming work will allow us to compute the local heights away from $p$, so we can use quadratic Chabauty here, too.},177,205,206,213,221,255,266,287,299,330 \}
\end{equation}
and the seven levels
\begin{equation}
\label{ggt2_hypell_levels}
    \{ 136, 171, 207, 252, 315, 176, 279 \}
\end{equation}
for which $X_0(N)^*$ is hyperelliptic of genus $g \geq 3$. For the square-free levels we also compute the $j$-invariants of the $\Q$-curves associated to non-cuspidal rational points. The methods of \cite{BGX21}, quadratic and bielliptic quadratic Chabauty, and recent work of Arul and M{\"u}ller~\cite{ArulMueller} together determine $X_0(N)^*(\Q)$ for all of the $18$ non-squarefree levels $N$ as explained in Section~\ref{subsection_nonsqfreelevels}. Thus, this completes the determination of rational points on all hyperelliptic curves $X_0(N)^*$. Together with the results in~\cite{BGX21,BBBLMTV21,BDMTV21,ArulMueller,AABCCKW,BianchiPadurariu}, we have the following theorem.

\begin{namedtheorem}[Main]
\label{maintheorem}
Let $N$ be such that $X_0(N)^*$ is hyperelliptic. Then the set $X_0(N)^*(\Q)$ consists only of points of `small' height.

Suppose $N$ is a square-free positive integer such that $X_0(N)^*$ is hyperelliptic. Then Table~\ref{tab:jinv} classifies each rational point on $X_0(N)^*$ as cusp, CM or exceptional. If $X_0(N)^*$ has no exceptional rational points, then $N \in \{ 67,107,146,167,205,213,390 \}$. For each of the remaining $32$ levels $N \in \{ 73$, $85$, $93$, $103$, $106$, $115$, $122$, $129$, $133$, $134$, $154$, $158$, $161$, $165$, $166$, $170$, $177$, $186$, $191$, $206$, $209$, $215$, $221$, $230$, $255$, $266$, $285$, $286$, $287$, $299$, $330$, $357 \}$, there is at least one exceptional rational point.
\end{namedtheorem}

\begin{remark}
\label{jinvtableremark}
For each of the $39$ square-free levels when $X_0(N)^*$ is hyperelliptic, Table~\ref{tab:jinv} in Appendix~\ref{sec:appendix} gives more refined information about the rational points on $X_0(N)^*$. It gives the $j$-invariants of the $\Q$-curves associated to all non-cuspidal rational points, and the CM discriminants for all the CM points.
\end{remark}

\begin{remark}
The boundedness conjecture of Elkies~\cite{Elkies2004} predicts that there are only finitely many positive integers $N$ such that $X_0(N)^*$ contains an exceptional rational point. Although the main theorem of this paper says that exceptional rational points exist on most of the hyperelliptic curves $X_0(N)^*$ of squarefree level $N$, we observe from Table~\ref{tab:jinv} that most of these exceptional points arise as the image of a cusp or CM point under the hyperelliptic involution. The only curves that have an exceptional rational point not arising in this way are $X_0(129)^*$ and $X_0(286)^*$. Furthermore, the curve $X_0(129)^*$ has extra automorphisms which in fact explain all the exceptional rational points on this curve.
\end{remark}

\begin{remark}
See Table~\ref{tab:curvessolvedpreviously} for the list of curves $X_0(N)^*$ whose rational points were determined in previous works. See Table~\ref{tab:curvessolvedhere} for the exact list of curves $X_0(N)^*$ whose rational points are determined in this paper.
\end{remark}

We did not make major improvements to the existing algorithms. The main challenge in our computations was to adapt the existing algorithms to our situations and to carefully choose suitable parameters, e.g., the primes in the various Chabauty methods and in the Mordell-Weil sieve, as well as the modulus in the Mordell-Weil sieve.

In Sections~\ref{ssec:Chabauty-Coleman} to~\ref{ssec:excluding one bad Qp point}, we give an overview of the ingredients including the various Chabauty methods and the Mordell-Weil sieve. We indicate in~\ref{ssec:quotients} how to compute the covering morphisms
\[
    X_0(N) \to X_0(N)/W'(N) \to X_0(N)^*
\]
for all subgroups $W'(N)$ of $W(N)$, which allows one to compute the $\Q$-rational points on the coverings, too. In Section~\ref{ssec:runtimes}, we give the runtimes for intensive quadratic Chabauty computations. In Sections~\ref{subsection_sqfreelevels} and~\ref{subsection_nonsqfreelevels}, we give details on the computations for the genus $2$ curves $X_0(N)^*$, with $N$ squarefree and non-squarefree respectively. Section~\ref{section_genusgt2} contains details of the computation for the curves $X_0(N)^*$ of genus $\geq 3$. Appendix~\ref{sec:appendix} contains a table giving refined information about the rational points on $X_0(N)^*$ for squarefree levels, as mentioned in the main theorem and Remark~\ref{jinvtableremark}.

\textbf{Acknowledgements.} We thank Jennifer Balakrishnan, Barinder Banwait, Lea Beneish, Abbey Bourdon, Netan Dogra, Sachi Hashimoto, Daniel Hast, Steffen Müller, Filip Najman, Michael Stoll, John Voight, and Boya Wen for helpful discussions. We are also grateful for the five anonymous referees' comments, which helped to significantly improve the article.

\section{Methods and algorithms used} \label{sec:algorithms}
Most of our computations rely on extensions of Chabauty's method and additional arguments (e.g., a variation of Mordell-Weil sieve). In this section, we give an overview of the methods used.

A curve over a field is called \emph{nice} if it is smooth, projective, and geometrically integral. Let $C/\Q$ be a nice curve of genus $g \geq 2$ with Jacobian $J$ whose Mordell-Weil group $J(\Q)$ has rank $r$.  For all flavors of Chabauty, assume that we have chosen a base point $x_0 \in C(\Q)$ to embed $C$ into $J$. In all our cases, there is such an $x_0$ because there is at least one $\Q$-rational cusp. If there were not one, then one could try to use the Mordell-Weil sieve to prove that $C(\Q)$ is empty.

To compute the rank of the Jacobian of a modular curve (over $\Q$), we use Magma's functionality to decompose the Jacobian into simple factors. For a simple modular abelian variety $A/\Q$, we can determine whether its rank $r \colonequals \rk A(\Q)$ equals $0$ or its dimension using~\cite{GrossZagier1986} and~\cite{KolyvaginLogachev} in all our cases.
\begin{theorem}[Gross--Zagier--Kolyvagin--Logach\"ev]
Let $A$ be a simple modular abelian variety over $\Q$ of dimension $g$ with algebraic rank $r \colonequals \rk A(\Q)$ and analytic rank $r_\mathrm{an} \colonequals \ord_{s=1}L(A,s) \in \{0,g\}$. Then $r = r_\mathrm{an}$, i.e., the BSD rank conjecture holds for $A$.
\end{theorem}
It is sufficient to check that $\ord_{s=1}L(f,s) \in \{0,1\}$ for one of the newforms associated with $A$~\cite[Corollary~V.1.3]{GrossZagier1986}. One can either read off the analytic rank of $f$ from the LMFDB~\cite{lmfdb} or use modular symbols to decide whether $L(f,1) = 0$; to prove that $L'(f,1) \neq 0$, one can use Dokchitser's code~\cite{DokchitserLvalues} for computing special $L$-values.

We used Magma~\cite{Magma} for our Chabauty computations of all flavors and Sage~\cite{sagemath} to access Liu's \texttt{genus2reduction} from \cite{Liu,Liu2}. The code to reproduce our computations can be found in the GitHub repository~\cite{ACKP}. The log files were produced with Magma~2.26-10.

First we give an overview on the $27$ levels $N$ such that $X_0(N)^*(\Q)$ was determined prior (or concurrently) to this paper. N.~B.~the level $161$ appears twice.
{ \small
\begin{table}[H]
\centering
\captionsetup{justification=centering,margin=0.1cm}
\begin{tabular}{ccc}
\toprule
Method & Levels $N$ & Reference \\ \midrule
Elliptic curve Chabauty & \begin{tabular}{@{}c@{}}85, 93, 106, 115, 122, 129, \\ 154, 158, 161, 165, 170, 186, 209,\\ 215, 230, 285, 286, 357, 390\end{tabular} & \cite{BGX21}\\ \midrule
 & 67, 73, 103 & \cite{BBBLMTV21} \\
Quadratic Chabauty & 107, 161, 167, 191 &  \cite{BDMTV21} \\
& 125 & \cite{ArulMueller} \\ \midrule
Bielliptic quadratic Chabauty & 166 & \cite{BianchiPadurariu} \\
\bottomrule
\end{tabular}
\caption{Levels $N$ and methods previously used to determine $X_0(N)^*(\Q)$}
\label{tab:curvessolvedpreviously}
\end{table}
}

Now we give an overview on the $37$ levels $N$ discussed in this paper and methods we used to determine $X_0(N)^*(\Q)$.
\begin{table}[H]
\centering
\captionsetup{justification=centering,margin=0.4cm}
\begin{tabular}{cc}
\toprule
Method & Levels $N$ \\ \midrule
Classical Chabauty & \begin{tabular}{@{}c@{}} 88, 104, 112, 116, 117, 121, \\  135, 136, 153, 168, 171, 176, 180, \\ 184, 198, 204, 276, 279, 284, 380 \end{tabular} \\ \midrule
Exceptional isomorphisms & 134, 146, 206, 266 \\ \midrule
Elliptic curve quotient & 207, 252, 315 \\ \midrule
Elliptic curve Chabauty & 147, 255, 330 \\ \midrule
Quadratic Chabauty & 133, 177, 205, 213, 221, 287, 299 \\
\bottomrule
\end{tabular}
\caption{Levels $N$ and methods we applied to determine $X_0(N)^*(\Q)$}
\label{tab:curvessolvedhere}
\end{table}
The code for all our computations is adapted from previous papers cited in the corresponding subsections below.

\subsection{The classical Chabauty--Coleman method} \label{ssec:Chabauty-Coleman}

In cases where $r < g$, we often use the classical method of Chabauty--Coleman. A good overview article on this method is~\cite{McCallumPoonen} and details on a refinement of this method can be found in~\cite{stoll_2006}.

Choose a prime $p$ of good reduction for $C$. The set $C(\Q)$ is contained in the set
\[
    C(\Q_p)_1 \colonequals \Bigg\{x \in C(\Q_p) : \int_{x_0}^x\omega = 0\Bigg\} \subseteq C(\Q_p)
\]
given as the vanishing locus of specific power series given by Coleman integration of a certain differential $\omega \in \H^0(C_{\Q_p}, \Omega^1) \iso \H^0(J_{\Q_p},\Omega^1)$ with the property that $\int_{[D]}\omega = 0$ for all $[D] \in J(\Q)$. For the construction and properties of the Coleman integral on curves, see \cite{ColemanIntegrals}, and \cite{BBK} for how to compute them on hyperelliptic curves. By linearity of the Coleman integral, it is enough to check the vanishing of $\int_{[D]}\omega$ for all $[D]$ from a generating set of a finite index subgroup of $J(\Q)$. The \emph{Chabauty condition} $r < g$ is used to construct a non-zero such differential $\omega$, which implies that $C(\Q_p)_1$ and hence $C(\Q)$ is finite. (In fact, it is enough that the $\Z_p$-rank of the $p$-adic closure of $J(\Q)$ in $J(\Q_p)$ is less than $g$.) This is because the integral will be non-zero and a non-zero power series has only finitely many zeros on a $1$-dimensional $p$-adic manifold.

It is crucial that one can compute the $p$-adic integral as a locally analytic function and hence explicitly approximate the Chabauty set $C(\Q_p)_1$. If all known $\Q$-points coincide with the $p$-adic approximations, then one has determined $C(\Q)$. If $g = 2$ and $r = 1$, then we use a combination of Chabauty--Coleman and the Mordell-Weil sieve, which is included in Magma. For $g \geq 3$, we use the code of Balakrishnan and Tuitman in~\cite{BalTui} to perform explicit Coleman integration for smooth curves. It enables us to determine the points in $C(\Q_p)_1$ precisely enough to check that this set consists exactly of the known rational points.

The code assumes that the degree-$0$ divisors given by the difference of the known $\Q$-points of $C$ generate a finite index subgroup of $J(\Q)$. For the levels $N = 171, 176$, the Jacobian of the curve $X_0(N)^*$ has rank $1$, so we simply confirm that one divisor has infinite order. For level $N = 279$, we check this by looking for at the reduction morphism to $J(\F_\ell)$ for small primes $\ell$ of good reduction to prove that there are no non-trivial relations between the divisor classes. For details, confer~\cite[Lemma~4]{Stoll6}.

\subsection{Elliptic curve Chabauty} \label{ssec:EllipticCurveChabauty}

When the Chabauty condition is violated, elliptic curve Chabauty is a method that sometimes succeeds in computing $C(\Q)$ using a combination of descent (also called going-up in \cite[\S{5.5}]{BjornExpo}) as described in~\cite[\S{1}]{StollDescent} and an analogue of the Chabauty--Coleman method over number fields. This method is described in~\cite{BruinEC,BruinFlynn}, and in~\cite{BGX21}.
We give a short summary in the way we use it.

Assume that $\pi: C \to \bbP^1$ is a finite morphism of degree $2$ and that $C$ is a hyperelliptic curve of genus $g = 2$. Let $\alpha: C' \to C$ be a twist of the pullback of the multiplication-by-$2$ morphism on $J$. (For the definition of a twist of a finite \'etale Galois covering, see~\cite[\S{4}]{StollDescent}.) Note that these are finite étale coverings with geometric Galois group $(\Z/2\Z)^4$. One has $C(\Q) = \bigcup_\alpha \alpha\left(C'(\Q)\right)$, where $\alpha$ runs over all such twists (``descent''). Hence we only need to determine $C'(\Q)$ for the finitely many curves $C'$ such that $C'$ has points locally everywhere. These can be computed from the fake $2$-Selmer set as described in~\cite{BruinStoll2009} and in Magma using \texttt{TwoCoverDescent}. Note that the genus of the coverings will grow, so the Chabauty condition is more likely to be satisfied at the expense of increasing the degree.

If $C'(\Q) = \emptyset$, we can use the Mordell-Weil sieve to prove $C'(\Q) = \emptyset$. If $C'(\Q) \neq \emptyset$, we can compute $C'(\Q)$ using the Chabauty method for an elliptic curve quotient over a finite extension of $\Q$. Since $C' \to C$ is Galois with Galois group $(\Z/2\Z)^{4}$ (as is $[2]: J \to J$), $\Jac(C')$ is isogenous to $J \times A$ with the Prym variety $A \colonequals \prod_{i=1}^{15} E_i$ with elliptic curves $E_i$ over $\ol\Q$ (one has $15 = 2^{2g} - 1$). These come equipped with a morphism $\pi': E_i \to \bbP^1$. Since $g=2$, every $\Z/2$ cover of $C_{\ol{\Q}}$ will have an elliptic Prym variety, and $A = \Res_{K/\Q}E \times_{\Q} \ol\Q$ for some elliptic curve $E$ over the \'etale algebra $K$ of $J[2] \setminus \{0\}$.  Now one can apply the classical Chabauty--Coleman method to $E$ if $\rk A(\Q) = \rk E(K) < [K:\Q] = \dim A$ (note that the dimension of a Jacobian of a curve equals the genus of that curve), together with the information that the images of the $\Q$-rational points on $C'$ on $E$ have image in $\bbP^1(\Q)$ under $\pi'$.

We have replaced an abelian variety over $\Q$ of dimension $d$ by an elliptic curve over a (possibly) larger field of degree $d$. One usually suspects this method to work if the degree of $K$ is large because the ranks of elliptic curves are expected to be small, assuming that one can compute them. However, elliptic curve Chabauty often fails for algorithmic reasons. Namely, if the \'etale algebra $K$ has too large a degree, to perform the $2$-descent to compute $A(\Q)$, one would need to compute the class group of a number field of large degree. Such a computation is infeasible, often even assuming GRH. In these cases we can use the quadratic Chabauty method, see Section~\ref{ssec:QC}.

We use this method to determine $X_0(N)^*(\Q)$ for $N = 147$, $255$ and $330$. In Section~\ref{subsection_nonsqfreelevels} we work out the case $N = 147$ in more detail. Our computations were performed using code written by \cite{BGX21}.

\subsection{Quadratic Chabauty} \label{ssec:QC}

If the $p$-adic closure of the image of $J(\Q) \subset J(\Q_p)$ under $\log$ has rank $g$, the classical Chabauty condition is violated. In this case, we can try to use quadratic Chabauty. %, for a prime $p$ of good reduction, the abelian logarithm induces a homomorphism $\log \colon J(\Q_p) \rightarrow \H^0(C_{\Q_p}, \Omega^1)^{\vee}$.

If the Néron-Severi rank $\rho(J)$ of $J$ is larger than $1$, then there exists a non-trivial $Z \in \Ker(\NS(J) \rightarrow \NS(C))$ inducing a correspondence on $C \times_\Q C$. Balakrishnan, Dogra, M\"uller, Tuitman and Vonk in~\cite{BD1,BDMTV} explain how to attach to any such $Z$ a locally analytic \emph{quadratic Chabauty function}
\[
    \rho_Z : C(\Q_p)\rightarrow \Q_p
\]
as follows: using Nekovář's theory of $p$-adic heights \cite{Nekovar}, one can construct a global $p$-adic height which decomposes as a sum of local height functions; for the latter, see \cite{ColemanGross}. The quadratic Chabauty function $\rho_Z$ is defined as the difference between the global $p$-adic height and the local height for the chosen prime $p$. Their difference equals the sum of the local heights away from $p$. % Even though we do not go into details here, we note that the computation of the global height pairing is easier when $C$ has many rational points.

The local $p$-adic heights $h_v$ with $v \neq p$ take only finitely many values on $C(\Q_v)$, see \cite{KimTamagawa}. Therefore, there exists a finite set $\Upsilon \subset \Q_p$, such that $\rho_Z(x) \in \Upsilon$ for any $x \in C(\Q)$. As described in \cite{BettsDogra}, this set is often computable a priori by intersection theory.

Since $\rho_Z$ has Zariski-dense image on every residue disk and is given by a convergent power series, this implies that the set of $\Q_p$-rational points of $C$ having values in $\Upsilon$ under $\rho_Z$ is finite. As both $C(\Q_p)_2$ and $C(\Q)$ are contained in that set, $C(\Q)$ is finite as well. Since both $\rho_Z$ and $\Upsilon$ can be explicitly computed by \cite{BDMTV,BDMTV21,BBBLMTV21}, this makes the provable determination of $C(\Q_p)_2$ possible if the \emph{quadratic Chabauty condition} $r < g + \rho(J) - 1$ is satisfied. This works at least under the assumption that $\rho_Z$ has no repeated roots; otherwise, we get a finite superset. The issue is distinguishing between two $p$-adic roots which are close together and one repeated root. If we do not know whether a root is repeated or not, then there is no terminating process for determining the roots. Looking at the $p$-adic approximations to a high enough precision, it is possible to prove that there are no multiple roots if there are not, but it is not possible to prove that there is a multiple root. As we work up to finite precision, we can say that we either have one point or two which are close to each other.

The quadratic Chabauty condition holds for our curves since $r = g$ and $\rho(J) = \dim{J}$ because $J$ has real multiplication: There is an isomorphism $\NS(J) \otimes \Q \isoto \End(J)^\dagger \otimes \Q$ with $(-)^\dagger$ the Rosati involution associated to the canonical principal polarization of $J$~\cite[Proposition~17.2]{MilneAbelianVarieties}. Since $J$ has real multiplication, the Rosati involution acts trivially on $\End(J)$ and $\End(J)$ is a number field of degree equal to $\dim{J}$.

The quadratic Chabauty algorithm in~\cite{BDMTV} is specific to modular curves only when determining the nontrivial class $Z$, which is computed by the Hecke operator $T_p$ (determined by the Eichler-Shimura relation). The input consists of a plane affine patch of a modular curve $C/\Q$ (satisfying certain technical properties) and a prime $p$ of good reduction such that the Hecke operator $T_p$ generates $\End(J) \otimes_\Z \Q$. To choose the primes $p$, we tried the first primes satisfying this condition and picked those with which we could compute all the points using the Mordell-Weil sieve. The algorithm outputs a finite set of $\Q_p$-points to some precision containing $C(\Q_p)_2$. 

In all our examples where we apply quadratic Chabauty, the set $\Upsilon$ equals $\{0\}$. This is because the stable reduction type at all primes of bad reduction computed with Liu's \texttt{genus2reduction} is the same as the one in~\cite[Example~5.19]{BDMTV21} with $v = 23$. Namely, it is $[\mathrm{I}_{1\text{-}1\text{-}0}]$ in the notation of~\cite{NamikawaUeno}, see p.~179 there for this case: The stable model has special fiber a curve of genus $0$ with exactly two double points and the given equation defines a regular semistable model over $\Z_v$ for all bad primes $v$. We prove regularity of the given hyperelliptic model over $\Z$ in the singular points of the special fiber. We do this by checking that the $v$-valuation of the constant term of the defining equation with the coordinates of the singular point plugged in is $1$.

The semistability can be checked by invoking Magma's \texttt{IsNode} for all singular points of the special fiber. Hence all $\Q_v$-points lie on exactly one irreducible component of a minimal regular model over $\Z_v$, so by~\cite[Theorem~3.2]{BDMTV21}, the height contribution at $v$ is $0$. Note that all levels which we apply this to are odd, see Table~\ref{tab:genus2}, and the function \texttt{genus2reduction} computes regular semistable model at all primes $v\neq 2$. This function was developed in \cite{Liu} and \cite{Liu2}, and can be called from Sage.

Before running the quadratic Chabauty algorithm on a hyperelliptic curve, we apply an automorphism such that there are no $\Q_p$-rational points at infinity. In all our examples we could find such an automorphism by trying matrices with small entries in $\GL_2(\Z)$.

If $C: y^2 = f(x)$ is hyperelliptic, the \emph{bad} $\Q_p$-rational points for $p > 2$ are exactly those in a residue disc of a Weierstraß point over $\F_p$. We usually choose our quadratic Chabauty prime(s) $p$ such that there is no bad $\Q_p$-rational point. However, for the levels $N = 213, 287, 299$, we had to choose the quadratic Chabauty prime(s) $p$ such that there is exactly one bad $\Q_p$-rational point because they had only $4$ known $\Q$-rational points. This would not be enough to solve for the global height pairing, so we had to compute local heights at $p$ manually before running \texttt{QCModAffine}. This can currently be done only for odd degree hyperelliptic curves, which have a Weierstraß point at infinity, and these are bad. Forthcoming work of Gajovi\'c could be used to do the computations for even degree hyperelliptic curves. However, by the argument described in Section~\ref{ssec:excluding one bad Qp point}, we could conclude that this single bad $\Q_p$-rational point is not $\Q$-rational.

\subsection{Bielliptic quadratic Chabauty}

There is a variant of quadratic Chabauty useful for bielliptic curves described in~\cite{BD1,Bianchi2020} and being worked out with precision bounds and implemented by Fran\-ces\-ca Bianchi and the fourth author~\cite{BianchiPadurariu}. (The code to perform these computations will be published with their paper.) This method is applicable for $X_0(166)^*$ with bielliptic equation
\[
    y^2 = x^6 + 2x^4 + 17x^2 - 4.
\]
(Note that the polynomial on the right hand side factors into two irreducible factors over $\Q$, but elliptic curve Chabauty fails for this curve.) This curve is not further discussed in the current paper, but it will be part of a database of more than $400$ bielliptic genus $2$ curves from the LMFDB, whose rational points will be computed using this method. In this case, $r = g = 2$, so the quadratic Chabauty condition is fulfilled. The bielliptic quadratic Chabauty method is simpler as the algorithm makes use of the extra structure which bielliptic genus $2$ curves have, in particular the two quotient maps to elliptic curves of rank $1$. Most of the computations are then delegated to the level of the two elliptic curves. 

We have not attempted to apply the standard quadratic Chabauty method here, as the reduction type at the primes of bad reduction is more complicated, so the height contribution at the bad primes might be non-zero.

\subsection{The Mordell-Weil sieve}

The Mordell-Weil sieve was introduced in~\cite{BruinStoll} based on an idea by Scharaschkin \cite{Scharaschkin}. In the way we use it, it is also described in~\cite[\S{3.4}]{BDMTV21}, and we are using their code, with minor adaptations to our examples. Another reference with more details about the Mordell-Weil sieve is~\cite[\S{6.7}]{BBBLMTV21}.

For a finite set $S$ of primes and an integer $M > 1$, consider the commutative diagram
\[\begin{tikzcd}
    C(\Q) \ar[r] \ar[d] & J(\Q)/MJ(\Q) \ar[d,"\alpha"] \\
    \displaystyle\prod_{\ell \in S}C(\Q_\ell) \ar[r,"\beta"] & \displaystyle\prod_{\ell \in S}J(\Q_\ell)/MJ(\Q_\ell)
\end{tikzcd}\]
with the horizontal maps induced by a $\Q$-rational divisor class of degree~$1$. Assume that one can compute $J(\Q)$ in terms of generators and relations, e.g., if the genus of $C$ is $2$ and $\Sha(J/\Q)[2] = 0$ using Stoll's \texttt{MordellWeilGroup\allowbreak{}Genus2} in Magma (this function often succeeds also if $\Sha(J/\Q)[2] \neq 0$, and outputs if its result is correct). Then one can compute the maps $\alpha$ and $\beta$. If one wants to prove that $C(\Q)$ is empty, the aim is to find $S$ and $M$ such that the images of $\alpha$ and $\beta$ are disjoint. More generally, the sieve can be used to show that a given residue disc in $C(\Q_\ell)$ does not contain a $\Q$-rational point. Conjecturally, for any curve one should be able to provably determine all $\Q$-points on a curve~\cite[\S{2}, p.~4]{BruinStoll} and~\cite[Main Conjecture]{StollDescent}.

We map $C(\Q_\ell)$ and $J(\Q_\ell)$ to $C(\F_\ell)$ and $J(\F_\ell)$ via the reduction map for primes $\ell$ of good reduction (not necessarily prime to $M$). Its fibers are the \emph{residue discs}.

We now give more details on how we use the sieve. Our finite set $S$ consits of the \emph{Mordell-Weil sieve primes} $\ell$ (``MWS primes'' for short). It contains as a subset the \emph{quadratic Chabauty primes} $p$ (``QC primes'' for short in Table~\ref{tab:genus2}).

Regarding the QC primes: To make the explanation easier, assume that $J(\Q)$ is torsion-free. Let $p$ be a QC prime, i.e., we have computed a set of residue discs of $C(\Q_p)$ using the quadratic Chabauty algorithm in Section~\ref{ssec:QC}, which contain the $\Q$-points of $C$ and such that there is at most one $\Q$-point in each residue disc. Let $P \in C(\Q_p)$ be a point which we suspect is not a $\Q$-rational point from our quadratic Chabauty computation, and that $P$ is given to some finite precision $p^n$. If $P$ were $\Q$-rational, then there would exist integers $a_1, \ldots, a_r$ such that
\[
    [P - x_0] = a_1 P_1 + \cdots + a_r P_r \in J(\Q),
\]
where $[P-x_0] \in J(\Q)$ is the image of $P \in C(\Q)$ under the Abel-Jacobi map and $\{P_1, \ldots, P_r\}$ is a generating set of $J(\Q)$. The strategy is to obtain a contradiction working modulo $p^{n}$. To this end, we compute a tuple $(\Tilde{a}_1, \ldots, \Tilde{a}_r) \in \Z/p^{n}\Z$ satisfying $\Tilde{a}_i \equiv a_i \;(\bmod \; p^{n})$ using the linearity of Coleman integrals of holomorphic differentials (see \cite[equation (48)]{BBBLMTV21}). (Note that  \[\log_J\colon J(\Q_p)^1 \otimes_{\Z_p} \Q_p \isoto \H^0(J_{\Q_p},\Omega^1)^*, [D] \mapsto \Bigg(\omega \mapsto \int_{[D]}\omega\Bigg)\] given by Coleman integrals is an isomorphism. Here, $J(\Q_p)^1$ is the kernel of the reduction of $J(\Q_p)$ to $J(\F_p)$, and $J(\Q_p)^1$ is of finite index in $J(\Q_p)$.) This gives us a coset in $p^nJ(\Q)$. To prove that $P$ is not $\Q$-rational, it is enough to show that the image of the corresponding coset of $p^{n} J(\Q)$ under $\alpha$ is not hit by $\beta$. (Note that since we are taking the quotient of $J(\F_\ell)$ by $MJ(\F_\ell)$ for all $\ell \in S$, this will only give us information if there are MWS primes $\ell$ with $p \mid \#J(\F_\ell)/MJ(\F_\ell)$. Hence we take $M$ to be divisible by $p^n$.) Sometimes we combine this information for several QC primes $p$ with the Chinese remainder theorem.

Regarding the MWS primes: From our Chabauty computations we get residue discs in $C(\Q_p)$ that contain at most one $\Q$-point. (This is the second way in which we use the information from our quadratic Chabauty computation in the Mordell-Weil sieve.) For the MWS primes $\ell$ that are not QC primes, we compute $C(\F_\ell)$ and also get residue discs $C(\Q_\ell)$. %For those of them which do not contain a known $\Q$-point, we try to show that their image under $\beta_p\colon C(\F_\p) \to J(\F_p)/MJ(\F_p)$ is not hit by $\alpha_p\colon J(\Q)/MJ(\Q)$. (This is the second way in which we use the information from our quadratic Chabauty computation in the Mordell-Weil sieve.) For the MWS primes $\ell$ that are not QC primes, we compute $C(\F_\ell)$ and get residue discs $C(\Q_\ell)$, and perform the same check as for the QC primes.
Some of these residue discs (obtained for QC primes and other MWS primes) contain the known $\Q$-rational points. The other ones are called \emph{fake residue discs}, and we want to show that they do not contain a $\Q$-rational point. To this end, we show that their image under $\beta$ is not hit by $\alpha$. This can be decided because $\alpha$ and $\beta$ are computable. 

We first take $M$ to be the product of the QC primes $p$ (or $p^n$, but in all our computations $p$ is enough) and run the sieve. If this does prove that the known $\Q$-points of $C$ are all of them, we multiply $M$ by a small integer $M'$ such that the $\#(J(\F_\ell)/MJ(\F_\ell))$ have common factors for the MWS primes $\ell \in S$. We use the heuristics from~\cite[\S{3}]{BruinStoll} to compute appropriate MWS primes and $M'$.

\subsection{Excluding $\Q$-rational points in Weierstraß discs} \label{ssec:excluding one bad Qp point}

Assume that we know (from classical or quadratic Chabauty) that there is at most one $\Q$-rational point in the $\Q_p$-residue disc of a Weierstraß point over $\F_p$, $p$ odd and of good reduction. Also assume that the Weierstraß point is not $\Q$-rational. Then the number of $\Q$-rational points in the residue disc is even (because the orbits of an involution have cardinality either $1$ or $2$, and orbit size $1$ means that it is a fixed point), hence zero because it is $\leq 1$. (We thank Michael Stoll for this argument.)  

Note that a power series which describes a regular differential at a Weierstraß point is even, so it becomes odd after integrating. Hence the Chabauty bound for this disk will be even. As the Weierstraß point is always a zero of the integral, one can improve the bound by $1$ if one knows that the point is not rational.

We use this argument to prove that the remaining bad $\Q_p$-rational point is not $\Q$-rational for $X_0(N)^*$ and $N = 213, 287, 299$. Here we had to allow one bad $\Q_p$-rational point because there were not enough $\Q$-rational points to solve for the height pairing and to apply the algorithm for the computation of the local Coleman--Gross $p$-adic heights at $p$. This algorithm, due to Balakrishnan and Besser, requires an odd degree $\Q_p$-model. (Forthcoming work of Gajovi\'{c} will allow one to compute the local heights also for even degree models.)

\subsection{Computation of Atkin-Lehner quotients} \label{ssec:quotients}

If the maximal quotient $X_0(N)^*$ has finitely many $\Q$-rational points and we know them, we can compute the $\Q$-rational points on all coverings $X_0(N)/W'(N)$ for a subgroup $W'(N)$ of $W(N)$ if we know the quotient morphism. In principle, one could compute the automorphism group of $X_0(N)$ over $\Q$ (or over $\F_p$ and try to lift the automorphisms) using Magma's \texttt{AutomorphismGroup} and the quotient using \texttt{CurveQuotient}. However, both functions turn out not to be very efficient for our values of $N$. Therefore we do the following:

Assume that $X_0(N)$ is of genus $g \geq 3$ and not hyperelliptic. We compute the canonical model of $X_0(N)$ in $\bbP^{g-1}$ as an intersection of quadrics by finding enough relations between monomials of degree $2$ of an integral basis of $S_2(\Gamma_0(N))$ until we get an irreducible $1$-dimensional scheme. (This is a well-known strategy, presented with more details in \cite{Galbraith_1996,BoxModels}.) Magma can compute the Atkin-Lehner operators as linear involutions on this basis. We compute the common fixed space of $\bbP^{g-1}$ by all Atkin-Lehner operators in a chosen subgroup $W'(N)$ of $W(N)$ by simultaneously diagonalizing them, so they are diagonal with entries $\pm1$ on their diagonal. Call this new basis $(g_i(q))_{i=1}^g$.

To compute the quotient for a given subgroup $W'(N)$ of the Atkin-Lehner involutions, let $(x_i)_{i=1}^g$ be the coordinates of the canonical embedding corresponding to $(g_i)_{i=1}^g$. Let $w_q$ act on $x_i$ with the same sign as it acts on $g_i$. Let $M$ be the set of all degree $2$ monomials $m_{i,j} = x_ix_j$ fixed by all Atkin-Lehner operators in $W'(N)$ and such that there is a $w_q \in W'(N)$ that does not fix $g_i$ or $g_j$. Let $J$ be the set of all $i$ such that $w_qg_i = g_i$ for all $w_q \in W'(N)$.  Let $\phi: \bbP^{g-1} \to \bbP_{(2,\ldots,2,1,\ldots,1)}$, $$[x_1:\ldots:x_g] \mapsto [m(x_1, \dotsc, x_g) (m \in M) : x_j (j \in J)]$$ (image in the weighted projective space). Here the first coordinates with $m \in M$ are homogeneous of degree $2$, and the last coordinates with $j \in J$ are homogeneous of degree $1$. Then by invariant theory, the quotient $X_0(N)/W'(N)$ is the image of the composition $X_0(N) \inj \bbP^{g-1} \to \bbP_{(2,\ldots,2,1,\ldots,1)}$. (We thank John Voight for this idea.)

There is another method to get the quotient $X_0(N) \to X_0(N)/W'(N)$ if the quotient is hyperelliptic:
Let $g_1(q),g_2(q) \in \Z\llbracket q\rrbracket$ be (the $q$-expansions of) two independent cusp forms with integral coefficients fixed by $W'(N)$. Assume we know the genus $g'$ of the quotient. Then, as in~\cite{BGX21}, we compute a hyperelliptic equation $y^2 = f(x)$ of the quotient as a solution to the linear equation system given by comparing coefficients in the hyperelliptic equation with $x(q) = \frac{g_2(q)}{g_1(q)}$ and $y(q) = qx'(q)\frac{1}{g_1(q)}$. The morphism from the canonical model of $X_0(N)$ to this hyperelliptic curve with weighted homogeneous coordinates $[x:y:z]$ of degree $(1,g'+1,1)$ is given on the $x$-coordinate by the coordinate corresponding to $g_2$ and on the $z$-coordinate by the coordinate corresponding to $g_1$. The $y$-coordinate is computed (assuming $g' = 2$) by writing $y(q)g_1(q)^3$ as a homogeneous polynomial of degree $g'+1 = 3$ in the $(g_i(q))_i$.

Since there are so many modular coverings, we actually computed the points only for one of them as a proof of concept, namely $X_0(133) \to X_0(133)^*$. It turns out that the $8$ $\Q$-rational points have $4$ $\Q$-rational preimages. This agrees with Kenku's computation~\cite[Theorem~1]{Kenku} of $X_0(N)(\Q)$, implying that $X_0(133)(\Q)$ consists only of cusps. Our code in the file \texttt{Coverings/X0N\_coverings.m} on GitHub can be used to do these computations for the other coverings.

\subsection{Runtimes and memory usage} \label{ssec:runtimes}

We record in Table~\ref{tab:genus2runtimes} the runtimes and the memory usage on (one core of) an AMD\,\textsuperscript{\tiny\textregistered} Ryzen Pro 3700U @\,2.30\,GHz with 13.6\,GB of RAM using Magma V2.26-10. All other computations had negligible runtimes of less then 20 seconds. The biggest contributing factor to the memory usage is the number of cosets to be sieved out in the Mordell-Weil sieve step.

\begin{longtable}{rrr}
\toprule
    $N$ & runtime in seconds & RAM used in MB\\
\midrule
\endfirsthead
\toprule
    $N$ & runtime in seconds & RAM used in MB\\
\midrule
\endhead
\bottomrule
\caption{Runtimes and memory usage in the cases where we applied quadratic and elliptic curve Chabauty, confer Table~\ref{tab:genus2} for quadratic Chabauty.}
\label{tab:genus2runtimes}
\endlastfoot
\bottomrule
\endfoot
     133 & 110  & 130 \\
     177 & 22   & 130 \\
     205 & 546  & 406 \\
     213 & 16822& 225 \\
     221 & 80   & 130 \\
     287 & 1533 & 130 \\
     299 & 6669 & 3196 \\
\midrule
    147 & 37 & 85 \\
    255 & 47 & 96 \\
    330 & 25 & 85 \\
\end{longtable}

\section{Genus $2$ curves $X_0(N)^*$} \label{section_genus2}

As mentioned in the introduction, thanks to the works of~\cite{Hasegawa97,BGX21,BBBLMTV21,BDMTV21}, the only genus $2$ curves $X_0(N)^*$ that remain to be studied are the following:\\
\noindent\textbf{Case 1} (Square-free level):
$N \in \{ 133,134,146,166,177,$ $205,206,213,221,255,$ $266,287,299,330 \}$\\
\noindent\textbf{Case 2} (Non-squarefree level):
$N \in \{ 88, 104, 112, 116, 117, 121, 125, 135, 147,$ $153, 168, 180, 184, 198, 204, 276,$ $284, 380 \}$\\

\subsection{Square-free levels} \label{subsection_sqfreelevels}
The curves of levels $N = 134,146,206$ can be addressed by observing the exceptional isomorphisms
\[
    X_0(2p)^* \cong X_0(p)^+ \quad \text{for } p \in \{67,73,103\},
\]
and $X_0(p)^{+}(\mathbb Q)$ for these $p$ were computed in~\cite{BBBLMTV21} and~\cite{BDMTV21}.

The set $X_0(166)^*(\Q)$ is determined in \cite{BianchiPadurariu} using bielliptic quadratic Chabauty over $\Q_7$ and the extra points are sieved out using the primes $\{107, 4211 \}$. We compute $X_0(N)^*(\Q)$ for the levels $N = 255$ and $330$ using elliptic curve Chabauty (their hyperelliptic polynomials factor over $\Q$ as a product of a quadratic and a quartic polynomial).

We investigate the remaining cases using quadratic Chabauty. Note that $X_0(266)^* \cong X_0(133)^*,$ thus we are left with
\[
    N \in \{ 133,177,205,213,221,287,299 \}.
\]
In Table~\ref{tab:genus2}, we record the parameters used in the cases where we use quadratic Chabauty. For all levels $N$ in this table, we explain why there is exactly one cusp on $X_0(N)^*$. The number of cusps on $X_0(N)$ is \[\sum_{d \mid N}\phi(\gcd(d,N/d)) = \sum_{d \mid N} 1 = 2^{\omega(N)}.\] To get the curve $X_0(N)^*$, we quotient $X_0(N)$ by the full group $W(N)$ of order $2^{\omega(N)}$. The action of $W(N)$ on these cusps is transitive, so $X_0(N)^*$ has exactly one cusp.
We denote this single cusp by $\infty$.

The quadratic Chabauty method is based on height computations. To use this method, one has to compute the contributions at the primes of bad reductions. In the case of a semistable reduction type with just one component, the contribution is trivial by~\cite[Theorem~3.2]{BDMTV21}. Namely, the local heights away from $p$ factor through the reduction graph, which has only one vertex in this case. This is fortunately the case in all of the cases in which we apply quadratic Chabauty (for $N = 166, 255, 330$ this does not hold, and we approach these using bielliptic quadratic and elliptic curve Chabauty). Therefore, one can run a quadratic Chabauty computation analogous to the one carried out for $X_0(67)^+$.
We conclude that the $\Q$-rational points of $X_0(N)^*(\Q)$ for $N$ in~\eqref{geq2_sqfree_levels} are precisely those predicted by Table~1 of~\cite{BGX21}.

To classify the rational points and compute corresponding $j$-invariants, we use the method outlined in~\cite{BGX21}. Denote by $j(z)$ the usual $j$-function on $X_0(1)$. We find suitable generators of the function field of $X_0(N)^*$. Then we use them to express the symmetric functions of the set $\{j(dz) : 1 \leq d \mid N\}$, in terms of these generators. This is how we find the $j$-invariants of the $\Q$-curves associated to each non-cuspidal rational point on $X_0(N)^*$. Table~\ref{tab:jinv} in the appendix lists these $j$-invariants when they lie in $\Q$ or a quadratic number field, and the defining fields $\Q(j)$ otherwise. The minimal polynomials of all these $j$-invariants are given in the file \texttt{jinvs.log} in our GitHub repository~\cite{ACKP}. Table~\ref{tab:jinv} also classifies the points of $X_0(N)^*$ as CM points or exceptional points, and gives the discriminant $D$ of the order of the endomorphism ring when the point is CM.

\begin{longtable}{cccccc}
\toprule
     $N$ & QC primes & MWS primes & $M'$ & $\#X_0(N)^*(\Q)$ & $\#$non-CM\\
\midrule
\endfirsthead
\toprule
     $N$ & QC primes & MWS primes & $M'$ & $\#X_0(N)^*(\Q)$ & $\#$non-CM\\
\midrule
\endhead
\bottomrule
\caption{Data used to do the quadratic Chabauty computations and information on $X_0(N)^*(\Q)$. (MWS primes means the additional primes in the Mordell-Weil sieve compared to the QC primes.)}
\label{tab:genus2}
\endlastfoot
\bottomrule
\endfoot
     133 & 5, 59& 109, 131, 317, 509 & $1$ & 8 & 2\\
     177 & 5 & 19 & $1$ & 6 & 1\\
     205 & 17, 61, 71 & $\leq 18457$\footnote{43, 71, 179, 359, 439, 617, 661, 967, 1033, 1997, 2063, 2213, 2381, 2753, 3373, 9579, 15083, 18457} & $3$ & 6 & 0\\
     213 & 79\footnote{One can also use the QC primes $37$, $43$ with MWS primes $\leq 23819$ and $M' = 2^2 \cdot 3 \cdot 5 \cdot 29$.} & 59, 149, 211, 4177 & $1$ & 4 & 0\\
     221 & 29 & 3, 47 & $1$ & 6 & 1\\
     287 & 19, 29 & 3, 5 & $1$ & 4 & 1\\
     299 & 29, 37, 43 & $\leq 929$\footnote{7, 59, 89, 103, 137, 317, 443, 541, 787, 929} & $3 \cdot 5$ & 4 & 1\\
\end{longtable}

\subsection{Non-squarefree levels} \label{subsection_nonsqfreelevels}
Although only square-free levels $N$ are considered in~\cite{BGX21}, the methods of elliptic curve Chabauty and quadratic Chabauty work beyond this scope. Among the $18$ non-squarefree levels $N$ for which $X_0(N)^*$ has genus $2$, we have that $\text{rk } J_0(N)^*(\Q) \in \{0,1\}$ for all except for $N = 147$ and $N = 125$. Thus in those cases one can compute $X_0(N)^*(\Q)$ using the classical Chabauty method implemented in Magma. Recent work of Arul--M{\"u}ller~\cite{ArulMueller} verifies that the known small rational points on $X_0(125)^+= X_0(125)^*$ are in fact all of them.

For $N = 147$, we modified the code attached to~\cite{BGX21}, as it works to compute $X_0(147)^*(\Q)$ using elliptic curve Chabauty: An equation for $X_0(147)^*$ is given by
\[ y^2 = f(x)=(x^2 - x + 1)(x^4 + x^3 - 4x^2 - 3x + 9).\]
We let $K \colonequals \Q[t]/(t^4 + t^3 - 4t^2 - 3t + 9)$. 
Then over $K$ the polynomial $f(x)$ factorizes as
$$f(x) = L_1(x)L_2(x)Q_1(x)Q_2(x),$$
where the polynomials $L_i$ are linear, and the polynomials $Q_i$ are quadratic. Let $E$ be the elliptic curve $y^2 = Q_1(x)Q_2(x)$ over $K$.
If $(x_0,y_0)$ is in $X_0(147)^*(\Q)$, then $x_0$ is the $x$-coordinate of a point on $E(K)$ or on one of its $3$ twists we have to consider.  We computed these twists as in~\cite{BruinStoll2009} and the corresponding rational points in Magma using \texttt{TwoCoverDescent}.
We conclude that $X_0(147)^*(\Q)$ equals
\[ X_0(147)^*(\Q) \hspace{-0.15em} = \hspace{-0.15em} \left\{ \pm\infty, (0,\pm 3), (1,\pm 2), (-2,\pm 7), \left(\frac{3}{2},\pm \frac{21}{2^3}\right),\left(-\frac{5}{3},\pm\frac{154}{3^3}\right) \right\}\hspace{-0.2em}.\]
We handle the $N=255$ and $N=330$ cases similarly.

We now proceed with an example of how to compute the points on a modular covering, simpler than the method described in Section~\ref{ssec:quotients}:
A model of $X_0(147)^+$ can be computed in Magma, together with its automorphisms. Since there are only two, the non-trivial automorphism of $X_0(147)^+$ defined over $\Q$ must be the projection of $w_3  = w_{49}$, and hence we obtain an explicit map
\[
    X_0(147)^+ \rightarrow X_0(147)^*.
\]
As we know $X_0(147)^*(\Q)$, we can look at preimages of the points under the above map to find that $\#X_0(147)^+(\Q) = 4$. 

B.~Banwait and F.~Najman were the first to point out the need to compute $\Q$-rational points on $X_0(147)^+$ and $X_0(147)^*$. Their forthcoming work with the fourth author \cite{BNP} determines all the cyclic isogenies of elliptic curves over $\Q(\sqrt{213})$. Before considering this number field, they examined $\Q(\sqrt{-5})$, and for this choice of number field they had to consider the level $N = 147$.

For the following non-squarefree levels $N$ with $r \in \{0,1\}$ we are using classical Chabauty--Coleman: The rank~$0$ levels are $N = 104$, $117$, $180$, and $168$, and the rank~$1$ levels are
$N = 88$,
$116$,
$121$,
$153$,
$184$,
$198$,
$380$,
$112$,
$135$,
$204$,
$276$, and
$284$.
Arul and Müller~\cite{ArulMueller} determined the $\Q$-rational points on the curve $X_0(125)^+$ with rank $r = 2$ using quadratic Chabauty. 

\section{Hyperelliptic $X_0(N)^*$ of genus $\geq 3$} \label{section_genusgt2}

Recall from~\cref{hyperelliptic star quotients} that the only levels $N$ such that $X_0(N)^*$ is hyperelliptic of genus $\geq 3$ are $N = 136$, $171$, $207$, $252$, $315$, $176$, and $279$.

These $7$ levels are addressed with the following methods:

\begin{itemize}
    \item Chabauty--Coleman: $136$ (Stoll's bound~\cref{Stoll Chabauty}), $171, 176, 279$;
    \item a degree $2$ map to a rank $0$ elliptic curve quotient for $207, 252, 315$;
    \item a covering collection technique for $176$. (Note that this level appears twice.)
\end{itemize}

We provide some details below.

We use the following quantitative sharpening of Chabauty--Coleman by Stoll:

\begin{theorem}[{\cite[Corollary~6.7]{stoll_2006}}] \label[theorem]{Stoll Chabauty}
Let $C$ be a nice curve of genus $g \geq 2$. Let $r$ be the rank of its Jacobian over $\Q$. Let $p$ be a prime of good reduction for $C$. If $r < g$ and $p > 2r+2$, then
\[
    \#C(\Q) \leq \#C(\F_p) + 2r.
\]
\end{theorem}

In Table~\ref{tab:genusgt2} we list all curves for which we use the \texttt{effective\_chabauty} function from~\cite{BalTui}. For these curves, we first verify that differences of rational points on the curve generate a finite index subgroup of the Mordell-Weil group of the Jacobian. This is done in the file \texttt{Coleman/diff\_of\_ratpts.m} in our GitHub repository~\cite{ACKP}. We do not need the Mordell-Weil sieve step for these curves. These are the only curves for which we actually need to compute a Chabauty function because the bound in~\cref{Stoll Chabauty} is not sharp.

\begin{longtable}{ccccc}
\toprule
     $N$ & $g$ & $r$ & $p$ & $\#X_0(N)^*(\Q)$ \\
\midrule
\endfirsthead
\toprule
     $N$ & $g$ & $r$ & $p$ & $\#X_0(N)^*(\Q)$ \\
\midrule
\endhead
\bottomrule
\endfoot
\bottomrule
\caption{Data used to do the Chabauty--Coleman computations using the function \texttt{effective\_chabauty} from~\cite{BalTui}.}
\label{tab:genusgt2}
\endlastfoot
     $171$ & $3$ & $1$ & $5$ & $6$ \\
     $176$ & $4$ & $1$ & $3$ & $5$ \\
     $279$ & $5$ & $2$ & $5$ & $6$ \\
\end{longtable}
%There is exactly one cusp on $X_0(167)^*$ because $N = 167$ is square-free.

The other levels are approached as follows:

For $X_0(136)^*$,
we have $g = 3$ and $r = 0$. Using~\cref{Stoll Chabauty} with $p=3$ we conclude that $\#X_0(136)^*(\Q) = 4$.

For $N \in \{207,252,315 \}$, the curve $X_0(N)^*$ has an explicit degree $2$ map to an elliptic curve of rank $0$, so the rational points are computed as preimages.

For $N = 176$, in addition to the Chabauty--Coleman method presented in the table above, one can also use a covering collection technique, as described below:

We have that the Jacobian of
$$X_0(176)^*: y^2 = x(x^3 - 2x^2 + 2)(x^3 - 4x + 4)(x^3 + 2x^2 - 2)$$ 
has rank $1$ over $\Q$.
Let $f_1 = x(x^3 - 2x^2 + 2)$ and $f_2 = (x^3 - 4x + 4)(x^3 + 2x^2 - 2)$, whose resultant is $2^{12}$. We then have to consider the following covering collection
$$X_d : y_1^2 = df_1(x), \quad y_2^2 = df_2(x),\quad d \in \{ \pm 1, \pm 2\}.$$
Consider the curves 
$$Y_d : y^2 = df_2(x), \quad d \in \{ \pm 1, \pm 2\}.$$
The Jacobians of $Y_1$ and $Y_{-2}$ have rank~$1$ and those of $Y_{-1}$ and $Y_2$ have rank~$0$,
so we can use \texttt{Chabauty} and \texttt{Chabauty0}, respectively, in Magma to find all the rational points. Then we can find all the rational points on the covering collection $\{X_d\}$, and then show that $\#X_0(176)^*(\Q) = 5$.

\newpage
\begin{appendix}
\section{Classification of $\Q$-rational points on hyperelliptic $X_0(N)^*$ for $N$ squarefree} \label{sec:appendix}

For each squarefree level $N$ when the curve $X_0(N)^*$ is hyperelliptic, Table~\ref{tab:jinv} gives information about all rational points on the curve. A lot of this data also appears in~\cite{BGX21}. But as mentioned earlier, that article does not contain information for the three levels $166, 255$ and $330$. Since this work completes the determination of rational points on all hyperelliptic $X_0(N)^*$, we think that presenting information about all the hyperelliptic curves $X_0(N)^*$ of squarefree level $N$ in a single place is valuable. Table~\ref{tab:jinv} also corrects several typos in the table given in~\cite{BGX21}. When a $j$-invariant has degree greater than $2$, the table only gives the defining fields $\Q(j)$. For accessing the $j$-invariant itself, we refer the reader to the file \texttt{jinvs.log} in our GitHub repository~\cite{ACKP}.

The prime factorization of $N$ is given below each level. An entry ``new'' below the prime factorization signifies that this level has been done in this paper.

While we use several different models of the curves $X_0(N)^*$ for carrying out the computation of rational points, the data in the table corresponds to the model obtained in Magma as\\
\texttt{> X := SimplifiedModel(X0NQuotient(N, PrimeFactors(N)));}\\
(Note that $N$ is assumed to be squarefree.)
The model $X$ is given by an equation $y^2 = f(x)$ for a monic sextic polynomial $f(x)$. The unique cusp of $X_0(N)^*$ is one of the points at infinity on $X$ which we denote by $\infty$.

One can notice that several $j$-invariants appear in the table for multiple levels $N$. This is explained by the Heegner points on $X_0(N)$ associated to orders in imaginary quadratic fields. For example, suppose $\mathcal{O}$ is an order in an imaginary quadratic field $K$ with class number $1$. If all the primes dividing $N$ are split or ramified in $\mathcal{O}$, then there is an ideal $\mathfrak{a} \subseteq \mathcal{O}$ of norm $N$ and the elliptic curve $\mathbb{C}/\mathcal{O}$ admits a cyclic $N$-isogeny to $\mathbb{C}/\mathfrak{a}$. The theory of Complex Multiplication says that the $j$-invariant of $\mathbb{C}/\mathcal{O}$ is a rational integer. Thus, given such a $j$-invariant, it appears as the $j$-invariant of a CM point on $X_0(N)$ for all levels $N$ satisfying the condition mentioned above. Hence it also appears in our table as the $j$-invariant of a CM point on $X_0(N)^*$ for all these levels $N$. See~\cite[\S{3}]{Galbraith_1999} for an explanation of $\Q$-points on $X_0(N)^+$ with CM.

{\footnotesize
\begin{longtable}{ccccc}
\caption{Rational non-cuspidal points $X_0(N)^*(\Q) \setminus \{\infty\}$, and the $j$-invariants and CM discriminants $D$ of the associated $\Q$-curves.}
\label{tab:jinv}\\
\toprule
    $N$ & Point & CM & $D$ & $j$ or $\Q(j)$\\
\midrule
\endfirsthead
\toprule
    $N$ & Point & CM & $D$ & $j$ or $\Q(j)$\\
\midrule
\endhead
\bottomrule
\endfoot
\bottomrule
\endlastfoot
67 & $\infty'$ & $\text{yes}$ & $-11$ & $-2^{15}$\\
$67$ & $(-1, -7)$ & $\text{yes}$ & $-67$ & $-2^{15}~3^{3}~5^{3}~11^{3}$\\
 & $(-1, 7)$ & $\text{yes}$ & $-28$ & $3^{3}~5^{3}~17^{3}$\\
 & $(0, -3)$ & $\text{yes}$ & $-3$ & $0$\\
 & $(0, 3)$ & $\text{yes}$ & $-27$ & $-2^{15}~3^{1}~5^{3}$\\
 & $(1, -1)$ & $\text{yes}$ & $-7$ & $-3^{3}~5^{3}$\\
 & $(1, 1)$ & $\text{yes}$ & $-8$ & $2^{6}~5^{3}$\\
 & $(2, -1)$ & $\text{yes}$ & $-12$ & $2^{4}~3^{3}~5^{3}$\\
 & $(2, 1)$ & $\text{yes}$ & $-43$ & $-2^{18}~3^{3}~5^{3}$\\
\midrule
73 & $\infty'$ & $\text{yes}$ & $-12$ & $2^{4}~3^{3}~5^{3}$\\
$73$ & $(0, -1)$ & $\text{yes}$ & $-4$ & $2^{6}~3^{3}$\\
 & $(0, 1)$ & $\text{yes}$ & $-27$ & $-2^{15}~3^{1}~5^{3}$\\
 & $(1, -1)$ & $\text{yes}$ & $-8$ & $2^{6}~5^{3}$\\
 & $(1, 1)$ & $\text{yes}$ & $-19$ & $-2^{15}~3^{3}$\\
 & $(2, -3)$ & $\text{yes}$ & $-16$ & $2^{3}~3^{3}~11^{3}$\\
 & $(2, 3)$ & $\text{yes}$ & $-67$ & $-2^{15}~3^{3}~5^{3}~11^{3}$\\
 & $(\frac{3}{2}, \frac{-5}{8})$ & $\text{no}$ & $$ & $20 \left(2^{-26}{ 3 (-26670989 + 15471309 \sqrt{-127})}\right)^3$\\
 & $(\frac{3}{2}, \frac{5}{8})$ & $\text{yes}$ & $-3$ & $0$\\
\midrule
85 & $\infty'$ & $\text{yes}$ & $-19$ & $-2^{15}~3^{3}$\\
$5{\cdot}17$ & $(0, -5)$ & $\text{yes}$ & $-60$ & $1(3(470 + 213\sqrt{5}))^3(1 +
1\sqrt{5})/2$\\
 & $(0, 5)$ & $\text{yes}$ & $-35$ & $(16(-15 + 7\sqrt{5}))^3$\\
 & $(1, -2)$ & $\text{yes}$ & $-4$ & $2^{6}~3^{3}$\\
 & $(1, 2)$ & $\text{yes}$ & $-16$ & $2^{3}~3^{3}~11^{3}$\\
 & $(2, -5)$ & $\text{yes}$ & $-15$ & $1(3(-5 + 4\sqrt{5}))^3(3 - 
1\sqrt{5})/2$\\
 & $(2, 5)$ & $\text{yes}$ & $-115$ & $(48(-785 + 351\sqrt{5}))^3$\\
 & $(\frac{3}{2}, \frac{-17}{8})$ & $\text{no}$ & $$ & $\Q(\sqrt{17}, 
\sqrt{-95})$\\
 & $(\frac{3}{2}, \frac{17}{8})$ & $\text{yes}$ & $-51$ & $(48(5 - 
1\sqrt{17}))^3(-33 + 8\sqrt{17})$\\
 & $(\frac{-4}{3}, \frac{-425}{27})$ & $\text{yes}$ & $-595$ & $\Q(\sqrt{5}, 
\sqrt{17})$\\
 & $(\frac{-4}{3}, \frac{425}{27})$ & $\text{no}$ & $$ & $\Q(\sqrt{-4295}, 
\sqrt{-14603})$\\
\midrule
93 & $\infty'$ & $\text{yes}$ & $-12$ & $2^{4}~3^{3}~5^{3}$\\
$3{\cdot}31$ & $(-1, -3)$ & $\text{yes}$ & $-75$ & $(48(69 - 31\sqrt{5}))^3(0 + 
1\sqrt{5})$\\
 & $(-1, 3)$ & $\text{yes}$ & $-123$ & $(480(8 - 1\sqrt{41}))^3(-2049 + 
320\sqrt{41})$\\
 & $(0, -3)$ & $\text{yes}$ & $-24$ & $(12(5 + 2\sqrt{2}))^3(3 + 2\sqrt{2})$\\
 & $(0, 3)$ & $\text{yes}$ & $-60$ & $1(3(470 + 213\sqrt{5}))^3(1 + 
1\sqrt{5})/2$\\
 & $(1, -1)$ & $\text{yes, yes}$ & $-3, -27$ & $0, -2^{15}~3^{1}~5^{3}$\\
 & $(1, 1)$ & $\text{yes}$ & $-11$ & $-2^{15}$\\
 & $(2, -3)$ & $\text{yes}$ & $-15$ & $1(3(-5 + 4\sqrt{5}))^3(3 - 
1\sqrt{5})/2$\\
 & $(2, 3)$ & $\text{yes}$ & $-147$ & $(240(225 - 49\sqrt{21}))^3(-14 + 
3\sqrt{21})$\\
 & $(\frac{3}{2}, \frac{-9}{8})$ & $\text{no}$ & $$ & $\Q(\sqrt{-15}, 
\sqrt{-327})$\\
 & $(\frac{3}{2}, \frac{9}{8})$ & $\text{yes}$ & $-48$ & $(15(21 + 
13\sqrt{3}))^3(5 + 3\sqrt{3})$\\
 & $(\frac{1}{4}, \frac{-143}{64})$ & $\text{no}$ & $$ & $\Q(\sqrt{-23}, 
\sqrt{-143})$\\
 & $(\frac{1}{4}, \frac{143}{64})$ & $\text{yes}$ & $-3$ & $0$\\
\midrule
103 & $\infty'$ & $\text{yes}$ & $-67$ & $-2^{15}~3^{3}~5^{3}~11^{3}$\\
$103$ & $(0, -1)$ & $\text{yes}$ & $-27$ & $-2^{15}~3^{1}~5^{3}$\\
 & $(0, 1)$ & $\text{yes}$ & $-43$ & $-2^{18}~3^{3}~5^{3}$\\
 & $(1, -1)$ & $\text{yes}$ & $-11$ & $-2^{15}$\\
 & $(1, 1)$ & $\text{yes}$ & $-12$ & $2^{4}~3^{3}~5^{3}$\\
 & $(3, -19)$ & $\text{yes}$ & $-3$ & $0$\\
 & $(3, 19)$ & $\text{no}$ & $$ & $19(48(1623826405 + 30228849\sqrt{2885}))^3$\\
\midrule
106 & $\infty'$ & $\text{yes}$ & $-7$ & $-3^{3}~5^{3}$\\
$2{\cdot}53$ & $(-1, -4)$ & $\text{yes}$ & $-148$ & $(60(2837 + 
468\sqrt{37}))^3$\\
 & $(-1, 4)$ & $\text{yes}$ & $-36$ & $(4(21 + 20\sqrt{3}))^3(7 + 4\sqrt{3})$\\
 & $(0, -1)$ & $\text{yes, yes}$ & $-16, -4$ & $2^{3}~3^{3}~11^{3}, 
2^{6}~3^{3}$\\
 & $(0, 1)$ & $\text{yes}$ & $-40$ & $(6(65 + 27\sqrt{5}))^3$\\
 & $(1, -2)$ & $\text{yes}$ & $-52$ & $(30(31 + 9\sqrt{13}))^3$\\
 & $(1, 2)$ & $\text{yes}$ & $-24$ & $(12(5 + 2\sqrt{2}))^3(3 + 2\sqrt{2})$\\
 & $(2, -5)$ & $\text{yes}$ & $-4$ & $2^{6}~3^{3}$\\
 & $(2, 5)$ & $\text{yes}$ & $-100$ & $(6(2927 + 1323\sqrt{5}))^3$\\
 & $(\frac{1}{2}, \frac{-5}{8})$ & $\text{yes, yes}$ & $-28, -7$ & 
$3^{3}~5^{3}~17^{3}, -3^{3}~5^{3}$\\
 & $(\frac{1}{2}, \frac{5}{8})$ & $\text{no}$ & $$ & $\Q(\sqrt{33}, 
\sqrt{-2167})$\\
\midrule
107 & $\infty'$ & $\text{yes}$ & $-8$ & $2^{6}~5^{3}$\\
$107$ & $(0, -1)$ & $\text{yes}$ & $-43$ & $-2^{18}~3^{3}~5^{3}$\\
 & $(0, 1)$ & $\text{yes}$ & $-7$ & $-3^{3}~5^{3}$\\
 & $(2, -1)$ & $\text{yes}$ & $-28$ & $3^{3}~5^{3}~17^{3}$\\
 & $(2, 1)$ & $\text{yes}$ & $-67$ & $-2^{15}~3^{3}~5^{3}~11^{3}$\\
\midrule
115 & $\infty'$ & $\text{yes}$ & $-115$ & $(48(-785 + 351\sqrt{5}))^3$\\
$5{\cdot}23$ & $(1, -1)$ & $\text{yes}$ & $-11$ & $-2^{15}$\\
 & $(1, 1)$ & $\text{yes}$ & $-19$ & $-2^{15}~3^{3}$\\
 & $(2, -5)$ & $\text{yes}$ & $-15$ & $1(3(-5 + 4\sqrt{5}))^3(3 - 
1\sqrt{5})/2$\\
 & $(2, 5)$ & $\text{yes}$ & $-235$ & $(528(-8875 + 3969\sqrt{5}))^3$\\
 & $(\frac{1}{2}, \frac{-5}{8})$ & $\text{yes}$ & $-40$ & $(6(65 + 
27\sqrt{5}))^3$\\
 & $(\frac{1}{2}, \frac{5}{8})$ & $\text{no}$ & $$ & $\Q(\sqrt{65}, 
\sqrt{-9087})$\\
 & $(\frac{4}{3}, \frac{-35}{27})$ & $\text{no}$ & $$ & $\Q(\sqrt{10}, 
\sqrt{-9278})$\\
 & $(\frac{4}{3}, \frac{35}{27})$ & $\text{yes}$ & $-60$ & $1(3(470 + 
213\sqrt{5}))^3(1 + 1\sqrt{5})/2$\\
\midrule
122 & $\infty'$ & $\text{yes}$ & $-36$ & $(4(21 + 20\sqrt{3}))^3(7 + 
4\sqrt{3})$\\
$2{\cdot}61$ & $(-1, -4)$ & $\text{yes}$ & $-100$ & $(6(2927 + 
1323\sqrt{5}))^3$\\
 & $(-1, 4)$ & $\text{yes}$ & $-52$ & $(30(31 + 9\sqrt{13}))^3$\\
 & $(0, -1)$ & $\text{yes, yes}$ & $-16, -4$ & $2^{3}~3^{3}~11^{3}, 
2^{6}~3^{3}$\\
 & $(0, 1)$ & $\text{yes, yes}$ & $-12, -3$ & $2^{4}~3^{3}~5^{3}, 0$\\
 & $(1, -2)$ & $\text{yes}$ & $-20$ & $(2(25 + 13\sqrt{5}))^3$\\
 & $(1, 2)$ & $\text{yes}$ & $-88$ & $(60(155 + 108\sqrt{2}))^3$\\
 & $(\frac{3}{2}, \frac{-37}{8})$ & $\text{no}$ & $$ & $\Q(\sqrt{-15}, 
\sqrt{-951})$\\
 & $(\frac{3}{2}, \frac{37}{8})$ & $\text{yes}$ & $-232$ & $(30(140989 + 
26163\sqrt{29}))^3$\\
 & $(\frac{2}{3}, \frac{-37}{27})$ & $\text{no}$ & $$ & $\Q(\sqrt{-1598}, 
\sqrt{-1739})$\\
 & $(\frac{2}{3}, \frac{37}{27})$ & $\text{yes}$ & $-4$ & $2^{6}~3^{3}$\\
\midrule
129 & $\infty'$ & $\text{yes}$ & $-75$ & $(48(69 - 31\sqrt{5}))^3(0 + 
1\sqrt{5})$\\
$3{\cdot}43$ & $(-1, -3)$ & $\text{yes}$ & $-48$ & $(15(21 + 13\sqrt{3}))^3(5 + 
3\sqrt{3})$\\
 & $(-1, 3)$ & $\text{yes}$ & $-123$ & $(480(8 - 1\sqrt{41}))^3(-2049 + 
320\sqrt{41})$\\
 & $(0, -2)$ & $\text{yes}$ & $-8$ & $2^{6}~5^{3}$\\
 & $(0, 2)$ & $\text{yes}$ & $-147$ & $(240(-225 + 49\sqrt{21}))^3(14 - 
3\sqrt{21})$\\
 & $(1, -1)$ & $\text{yes}$ & $-12$ & $2^{4}~3^{3}~5^{3}$\\
 & $(1, 1)$ & $\text{yes, yes}$ & $-3, -27$ & $0, -2^{15}~3^{1}~5^{3}$\\
 & $(\frac{1}{2}, \frac{-3}{8})$ & $\text{no}$ & $$ & $\Q(\sqrt{57}, 
\sqrt{-687})$\\
 & $(\frac{1}{2}, \frac{3}{8})$ & $\text{yes}$ & $-51$ & $(48(-5 + 
1\sqrt{17}))^3(33 - 8\sqrt{17})$\\
 & $(\frac{-7}{5}, \frac{-383}{125})$ & $\text{no}$ & $$ & $\Q(\sqrt{-1059}, 
\sqrt{-135199})$\\
 & $(\frac{-7}{5}, \frac{383}{125})$ & $\text{yes}$ & $-3$ & $0$\\
 & $(\frac{7}{12}, \frac{-383}{1728})$ & $\text{no}$ & $$ & $\Q(\sqrt{85}, 
\sqrt{-347})$\\
 & $(\frac{7}{12}, \frac{383}{1728})$ & $\text{no}$ & $$ & $\Q(\sqrt{-7}, 
\sqrt{-444783})$\\
\midrule
133 & $\infty'$ & $\text{no}$ & $$ & $\Q(\sqrt{2}, \sqrt{69})$\\
$7{\cdot}19$ & $(0, -1)$ & $\text{yes}$ & $-19$ & $-2^{15}~3^{3}$\\
new & $(0, 1)$ & $\text{yes}$ & $-27$ & $-2^{15}~3^{1}~5^{3}$\\
 & $(1, -1)$ & $\text{yes}$ & $-12$ & $2^{4}~3^{3}~5^{3}$\\
 & $(1, 1)$ & $\text{yes}$ & $-91$ & $(48(-227 + 63\sqrt{13}))^3$\\
 & $(\frac{3}{5}, \frac{-83}{125})$ & $\text{no}$ & $$ & $\Q(\sqrt{-31}, 
\sqrt{-3651})$\\
 & $(\frac{3}{5}, \frac{83}{125})$ & $\text{yes}$ & $-3$ & $0$\\
\midrule
134 & $\infty'$ & $\text{yes}$ & $-52$ & $(30(31 + 9\sqrt{13}))^3$\\
$2{\cdot}67$ & $(-1, -3)$ & $\text{yes}$ & $-232$ & $(30(140989 + 
26163\sqrt{29}))^3$\\
new & $(-1, 3)$ & $\text{yes}$ & $-7$ & $-3^{3}~5^{3}$\\
 & $(0, -1)$ & $\text{yes, yes}$ & $-12, -3$ & $2^{4}~3^{3}~5^{3}, 0$\\
 & $(0, 1)$ & $\text{yes}$ & $-20$ & $(2(25 + 13\sqrt{5}))^3$\\
 & $(1, -1)$ & $\text{yes, yes}$ & $-28, -7$ & $3^{3}~5^{3}~17^{3}, 
-3^{3}~5^{3}$\\
 & $(1, 1)$ & $\text{yes}$ & $-8$ & $2^{6}~5^{3}$\\
 & $(\frac{-1}{2}, \frac{-7}{8})$ & $\text{yes}$ & $-72$ & $(20(49 + 
12\sqrt{6}))^3(49 + 20\sqrt{6})$\\
 & $(\frac{-1}{2}, \frac{7}{8})$ & $\text{no}$ & $$ & $\Q(\sqrt{113}, 
\sqrt{-1271})$\\
\midrule
146 & $\infty'$ & $\text{yes, yes}$ & $-12, -3$ & $2^{4}~3^{3}~5^{3}, 0$\\
$2{\cdot}73$ & $(-1, -1)$ & $\text{yes}$ & $-148$ & $(60(2837 + 
468\sqrt{37}))^3$\\
new & $(-1, 1)$ & $\text{yes}$ & $-36$ & $(4(21 + 20\sqrt{3}))^3(7 + 
4\sqrt{3})$\\
 & $(0, -1)$ & $\text{yes}$ & $-24$ & $(12(5 + 2\sqrt{2}))^3(3 + 2\sqrt{2})$\\
 & $(0, 1)$ & $\text{yes, yes}$ & $-16, -4$ & $2^{3}~3^{3}~11^{3}, 
2^{6}~3^{3}$\\
 & $(1, -3)$ & $\text{yes}$ & $-72$ & $(20(49 + 12\sqrt{6}))^3(49 + 
20\sqrt{6})$\\
 & $(1, 3)$ & $\text{yes}$ & $-8$ & $2^{6}~5^{3}$\\
 & $(2, -5)$ & $\text{yes}$ & $-4$ & $2^{6}~3^{3}$\\
 & $(2, 5)$ & $\text{yes}$ & $-100$ & $(6(2927 + 1323\sqrt{5}))^3$\\
\midrule
154 & $\infty'$ & $\text{yes}$ & $-40$ & $(6(65 + 27\sqrt{5}))^3$\\
$2{\cdot}7{\cdot}11$ & $(0, -2)$ & $\text{yes}$ & $-52$ & $(30(31 + 
9\sqrt{13}))^3$\\
 & $(0, 2)$ & $\text{yes}$ & $-24$ & $(12(5 + 2\sqrt{2}))^3(3 + 2\sqrt{2})$\\
 & $(1, -4)$ & $\text{yes, yes}$ & $-28, -7$ & $3^{3}~5^{3}~17^{3}, 
-3^{3}~5^{3}$\\
 & $(1, 4)$ & $\text{yes}$ & $-7$ & $-3^{3}~5^{3}$\\
 & $(2, 0)$ & $\text{yes}$ & $-84$ & $\Q(\sqrt{3}, \sqrt{7})$\\
 & $(\frac{-3}{2}, \frac{-77}{8})$ & $\text{yes}$ & $-1540$ & $\Q(\sqrt{5}, 
\sqrt{7}, \sqrt{11})$\\
 & $(\frac{-3}{2}, \frac{77}{8})$ & $\text{no}$ & $$ & $\Q(\sqrt{185}, 
\sqrt{2849}, \sqrt{-3367})$\\
 & $(\frac{-1}{3}, \frac{-56}{27})$ & $\text{yes, yes}$ & $-28, -112$ & 
$3^{3}~5^{3}~17^{3}, (15(2168 + 819\sqrt{7}))^3$\\
 & $(\frac{-1}{3}, \frac{56}{27})$ & $\text{no}$ & $$ & $\Q(\sqrt{7}, \sqrt{55},
\sqrt{-479})$\\
 & $(4, -22)$ & $\text{yes}$ & $-132$ & $\Q(\sqrt{3}, \sqrt{11})$\\
 & $(4, 22)$ & $\text{yes}$ & $-1848$ & $\Q(\sqrt{2}, \sqrt{33}, \sqrt{77})$\\
\midrule
158 & $\infty'$ & $\text{yes}$ & $-7$ & $-3^{3}~5^{3}$\\
$2{\cdot}79$ & $(0, -1)$ & $\text{yes}$ & $-24$ & $(12(5 + 2\sqrt{2}))^3(3 + 
2\sqrt{2})$\\
 & $(0, 1)$ & $\text{yes, yes}$ & $-12, -3$ & $2^{4}~3^{3}~5^{3}, 0$\\
 & $(2, -1)$ & $\text{yes}$ & $-148$ & $(60(2837 + 468\sqrt{37}))^3$\\
 & $(2, 1)$ & $\text{yes}$ & $-232$ & $(30(140989 + 26163\sqrt{29}))^3$\\
 & $(\frac{1}{2}, \frac{-1}{8})$ & $\text{yes, yes}$ & $-28, -7$ & 
$3^{3}~5^{3}~17^{3}, -3^{3}~5^{3}$\\
 & $(\frac{1}{2}, \frac{1}{8})$ & $\text{no}$ & $$ & $\Q(\sqrt{1169}, 
\sqrt{-1247})$\\
\midrule
161 & $\infty'$ & $\text{yes}$ & $-7$ & $-3^{3}~5^{3}$\\
$7{\cdot}23$ & $(-1, -7)$ & $\text{yes}$ & $-483$ & $\Q(\sqrt{69}, 
\sqrt{161})$\\
 & $(-1, 7)$ & $\text{yes}$ & $-91$ & $(48(-227 + 63\sqrt{13}))^3$\\
 & $(1, -1)$ & $\text{yes}$ & $-19$ & $-2^{15}~3^{3}$\\
 & $(1, 1)$ & $\text{yes}$ & $-115$ & $(48(-785 + 351\sqrt{5}))^3$\\
 & $(\frac{-1}{2}, \frac{-35}{8})$ & $\text{yes}$ & $-112$ & $(15(2168 + 
819\sqrt{7}))^3$\\
 & $(\frac{-1}{2}, \frac{35}{8})$ & $\text{no}$ & $$ & $\Q(\sqrt{-7}, 
\sqrt{32009})$\\
 & $(\frac{-1}{4}, \frac{-209}{64})$ & $\text{no}$ & $$ & $\Q(\sqrt{209}, 
\sqrt{-1140391})$\\
 & $(\frac{-1}{4}, \frac{209}{64})$ & $\text{yes}$ & $-28$ & 
$3^{3}~5^{3}~17^{3}$\\
\midrule
165 & $\infty'$ & $\text{yes}$ & $-11$ & $-2^{15}$\\
$3{\cdot}5{\cdot}11$ & $(0, -3)$ & $\text{yes}$ & $-51$ & $(48(-5 + 
1\sqrt{17}))^3(33 - 8\sqrt{17})$\\
 & $(0, 3)$ & $\text{yes}$ & $-195$ & $\Q(\sqrt{5}, \sqrt{13})$\\
 & $(1, 0)$ & $\text{yes}$ & $-24$ & $(12(5 + 2\sqrt{2}))^3(3 + 2\sqrt{2})$\\
 & $(2, -5)$ & $\text{yes}$ & $-35$ & $(16(-15 + 7\sqrt{5}))^3$\\
 & $(2, 5)$ & $\text{yes}$ & $-435$ & $\Q(\sqrt{5}, \sqrt{29})$\\
 & $(\frac{-1}{2}, \frac{-15}{8})$ & $\text{yes}$ & $-120$ & $\Q(\sqrt{2}, 
\sqrt{5})$\\
 & $(\frac{-1}{2}, \frac{15}{8})$ & $\text{no}$ & $$ & $\Q(\sqrt{-159}, 
\sqrt{265}, \sqrt{18497})$\\
 & $(-3, 0)$ & $\text{yes}$ & $-1155$ & $\Q(\sqrt{5}, \sqrt{33}, \sqrt{77})$\\
 & $(\frac{2}{3}, \frac{-55}{27})$ & $\text{no}$ & $$ & $\Q(\sqrt{-11}, 
\sqrt{-47}, \sqrt{661})$\\
 & $(\frac{2}{3}, \frac{55}{27})$ & $\text{yes, yes}$ & $-11, -99$ & $-2^{15}, 
(16(-77 + 15\sqrt{33}))^3(1057 - 184\sqrt{33})$\\
 & $(\frac{5}{2}, \frac{-99}{8})$ & $\text{no}$ & $$ & $\Q(\sqrt{-7}, \sqrt{33},
\sqrt{1441})$\\
 & $(\frac{5}{2}, \frac{99}{8})$ & $\text{yes}$ & $-1320$ & $\Q(\sqrt{5}, 
\sqrt{22}, \sqrt{33})$\\
\midrule
166 & $\infty'$ & $\text{yes}$ & $-8$ & $2^{6}~5^{3}$\\
$2{\cdot}83$ & $(0, -1)$ & $\text{yes}$ & $-20$ & $(2(25 + 13\sqrt{5}))^3$\\
new & $(0, 1)$ & $\text{yes}$ & $-52$ & $(30(31 + 9\sqrt{13}))^3$\\
 & $(1, -2)$ & $\text{yes}$ & $-24$ & $(12(5 + 2\sqrt{2}))^3(3 + 2\sqrt{2})$\\
 & $(1, 2)$ & $\text{yes}$ & $-88$ & $(60(155 + 108\sqrt{2}))^3$\\
 & $(-3, -10)$ & $\text{no}$ & $$ & $\Q(\sqrt{101}, \sqrt{2462})$\\
 & $(-3, 10)$ & $\text{yes}$ & $-232$ & $(30(140989 + 26163\sqrt{29}))^3$\\
 & $(\frac{-1}{3}, \frac{-10}{27})$ & $\text{no}$ & $$ & $\Q(\sqrt{-7523}, 
\sqrt{545254})$\\
 & $(\frac{-1}{3}, \frac{10}{27})$ & $\text{yes}$ & $-72$ & $(20(49 + 
12\sqrt{6}))^3(49 + 20\sqrt{6})$\\
\midrule
167 & $\infty'$ & $\text{yes}$ & $-43$ & $-2^{18}~3^{3}~5^{3}$\\
$167$ & $(-1, -1)$ & $\text{yes}$ & $-163$ & 
$-2^{18}~3^{3}~5^{3}~23^{3}~29^{3}$\\
 & $(-1, 1)$ & $\text{yes}$ & $-67$ & $-2^{15}~3^{3}~5^{3}~11^{3}$\\
\midrule
170 & $\infty'$ & $\text{yes}$ & $-36$ & $(4(21 + 20\sqrt{3}))^3(7 + 
4\sqrt{3})$\\
$2{\cdot}5{\cdot}17$ & $(-1, -2)$ & $\text{yes}$ & $-340$ & $\Q(\sqrt{5}, 
\sqrt{17})$\\
 & $(-1, 2)$ & $\text{yes, yes}$ & $-16, -4$ & $2^{3}~3^{3}~11^{3}, 
2^{6}~3^{3}$\\
 & $(0, -1)$ & $\text{yes}$ & $-15$ & $1(3(-5 + 4\sqrt{5}))^3(3 - 
1\sqrt{5})/2$\\
 & $(0, 1)$ & $\text{yes, yes}$ & $-4, -100$ & $2^{6}~3^{3}, (6(2927 + 
1323\sqrt{5}))^3$\\ \hdashline
 & $ (2, -5) $ & $\text{yes, yes}$ & $-15,-60$ &
\begin{tabular}{c}
$1(3(-5 + 4\sqrt{5}))^3(3 - 1\sqrt{5})/2$, \\
$1(3(470 + 213\sqrt{5}))^3(1 + 1\sqrt{5})/2$ \\
\end{tabular}\\ \hdashline
 & $(2, 5)$ & $\text{yes}$ & $-280$ & $\Q(\sqrt{2}, \sqrt{5})$\\
 & $(\frac{-1}{2}, \frac{-5}{8})$ & $\text{no}$ & $$ & $\Q(\sqrt{17}, \sqrt{65},
\sqrt{-247})$\\
 & $(\frac{-1}{2}, \frac{5}{8})$ & $\text{yes}$ & $-120$ & $\Q(\sqrt{2}, 
\sqrt{5})$\\
 & $(\frac{5}{3}, \frac{-38}{27})$ & $\text{yes}$ & $-4$ & $2^{6}~3^{3}$\\
 & $(\frac{5}{3}, \frac{38}{27})$ & $\text{no}$ & $$ & $\Q(\sqrt{-5}, \sqrt{19},
\sqrt{73})$\\
\midrule
177 & $\infty'$ & $\text{yes}$ & $-11$ & $-2^{15}$\\
$3{\cdot}59$ & $(0, -1)$ & $\text{yes}$ & $-8$ & $2^{6}~5^{3}$\\
new & $(0, 1)$ & $\text{yes}$ & $-24$ & $(12(5 + 2\sqrt{2}))^3(3 + 2\sqrt{2})$\\
 & $(\frac{3}{2}, \frac{-17}{8})$ & $\text{no}$ & $$ & $\Q(\sqrt{-23}, 
\sqrt{2881})$\\
 & $(\frac{3}{2}, \frac{17}{8})$ & $\text{yes}$ & $-267$ & $(240(-625 + 
53\sqrt{89}))^3(500001 - 53000\sqrt{89})$\\
\midrule
186 & $\infty'$ & $\text{yes, yes}$ & $-12, -3$ & $2^{4}~3^{3}~5^{3}, 0$\\
$2{\cdot}3{\cdot}31$ & $(-1, -3)$ & $\text{yes}$ & $-228$ & $\Q(\sqrt{3}, 
\sqrt{19})$\\
 & $(-1, 3)$ & $\text{yes}$ & $-84$ & $\Q(\sqrt{3}, \sqrt{7})$\\
 & $(0, -1)$ & $\text{yes}$ & $-24$ & $(12(5 + 2\sqrt{2}))^3(3 + 2\sqrt{2})$\\
 & $(0, 1)$ & $\text{yes}$ & $-15$ & $1(3(-5 + 4\sqrt{5}))^3(3 - 1\sqrt{5})/2$\\
 & $(1, -3)$ & $\text{yes}$ & $-120$ & $\Q(\sqrt{2}, \sqrt{5})$\\
 & $(1, 3)$ & $\text{yes}$ & $-168$ & $\Q(\sqrt{6}, \sqrt{14})$\\ \hdashline
 &$(2,-9)$&$\text{yes, yes}       $&$-15,-60$&
\begin{tabular}{c}
$1(3(-5 + 4\sqrt{5}))^3(3 - 1\sqrt{5})/2$, \\
$1(3(470 + 213\sqrt{5}))^3(1 + 1\sqrt{5})/2$ \\
\end{tabular}\\ \hdashline
 & $(2, 9)$ & $\text{yes}$ & $-708$ & $\Q(\sqrt{3}, \sqrt{59})$\\
 & $(\frac{-1}{2}, \frac{-3}{8})$ & $\text{yes, yes}$ & $-12, -48$ & 
$2^{4}~3^{3}~5^{3}, (15(21 + 13\sqrt{3}))^3(5 + 3\sqrt{3})$\\
 & $(\frac{-1}{2}, \frac{3}{8})$ & $\text{no}$ & $$ & $\Q(\sqrt{1257}, 
\sqrt{-2095}, \sqrt{24721})$\\
 & $(\frac{-4}{3}, \frac{-143}{27})$ & $\text{yes}$ & $-372$ & $\Q(\sqrt{3}, 
\sqrt{31})$\\
 & $(\frac{-4}{3}, \frac{143}{27})$ & $\text{no}$ & $$ & $\Q(\sqrt{37}, 
\sqrt{-143}, \sqrt{2077})$\\
\midrule
191 & $\infty'$ & $\text{yes}$ & $-19$ & $-2^{15}~3^{3}$\\
$191$ & $(0, -1)$ & $\text{yes}$ & $-7$ & $-3^{3}~5^{3}$\\
 & $(0, 1)$ & $\text{yes}$ & $-11$ & $-2^{15}$\\
 & $(2, -11)$ & $\text{yes}$ & $-28$ & $3^{3}~5^{3}~17^{3}$\\
 &$(2,11) $&$\text{no}$&& $j_0$\footnote[2]{$j_0 = \left(\frac{480}{191^2}(7725788647437 + 95942438\sqrt{2036079533}) \right)^3 \cdot (724537954586714121 + 16056976492100\sqrt{2036079533})$}\\
\midrule
205 & $\infty'$ & $\text{yes}$ & $-115$ & $(48(-785 + 351\sqrt{5}))^3$\\
$5{\cdot}41$ & $(0, -1)$ & $\text{yes}$ & $-40$ & $(6(65 + 27\sqrt{5}))^3$\\
new & $(0, 1)$ & $\text{yes}$ & $-16$ & $2^{3}~3^{3}~11^{3}$\\
 & $(-2, -7)$ & $\text{yes}$ & $-1435$ & $\Q(\sqrt{5}, \sqrt{41})$\\
 & $(-2, 7)$ & $\text{yes}$ & $-4$ & $2^{6}~3^{3}$\\
\midrule
206 & $\infty'$ & $\text{yes}$ & $-24$ & $(12(5 + 2\sqrt{2}))^3(3 + 
2\sqrt{2})$\\
$2{\cdot}103$ & $(-1, -1)$ & $\text{yes}$ & $-88$ & $(60(155 + 
108\sqrt{2}))^3$\\
new & $(-1, 1)$ & $\text{yes, yes}$ & $-12, -3$ & $2^{4}~3^{3}~5^{3}, 0$\\
 & $(0, -1)$ & $\text{yes}$ & $-20$ & $(2(25 + 13\sqrt{5}))^3$\\
 & $(0, 1)$ & $\text{yes}$ & $-40$ & $(6(65 + 27\sqrt{5}))^3$\\
 & $(\frac{1}{2}, \frac{-19}{8})$ & $\text{no}$ & $$ & $\Q(\sqrt{193}, 
\sqrt{-27119})$\\
 & $(\frac{1}{2}, \frac{19}{8})$ & $\text{yes}$ & $-148$ & $(60(2837 + 
468\sqrt{37}))^3$\\
\midrule
209 & $\infty'$ & $\text{yes}$ & $-8$ & $2^{6}~5^{3}$\\
$11{\cdot}19$ & $(0, -2)$ & $\text{yes}$ & $-88$ & $(60(155 + 108\sqrt{2}))^3$\\
 & $(0, 2)$ & $\text{yes}$ & $-19$ & $-2^{15}~3^{3}$\\
 & $(\frac{-1}{2}, \frac{-19}{8})$ & $\text{yes}$ & $-627$ & $\Q(\sqrt{33}, 
\sqrt{209})$\\
 & $(\frac{-1}{2}, \frac{19}{8})$ & $\text{no}$ & $$ & $\Q(\sqrt{-323551}, 
\sqrt{902537})$\\
\midrule
213 & $\infty'$ & $\text{yes}$ & $-51$ & $(48(5 - 1\sqrt{17}))^3(-33 + 
8\sqrt{17})$\\
$3{\cdot}71$ & $(1, -1)$ & $\text{yes}$ & $-11$ & $-2^{15}$\\
new & $(1, 1)$ & $\text{yes}$ & $-123$ & $(480(8 - 1\sqrt{41}))^3(-2049 + 
320\sqrt{41})$\\
\midrule
215 & $\infty'$ & $\text{no}$ & $$ & $\Q(\sqrt{2}, \sqrt{47645})$\\
$5{\cdot}43$ & $(1, -1)$ & $\text{yes}$ & $-19$ & $-2^{15}~3^{3}$\\
 & $(1, 1)$ & $\text{yes}$ & $-235$ & $(528(-8875 + 3969\sqrt{5}))^3$\\
 & $(2, -10)$ & $\text{yes}$ & $-115$ & $(48(-785 + 351\sqrt{5}))^3$\\
 & $(2, 10)$ & $\text{no}$ & $$ & $\Q(\sqrt{85}, \sqrt{11623937})$\\
\midrule
221 & $\infty'$ & $\text{yes}$ & $-16$ & $2^{3}~3^{3}~11^{3}$\\
$13{\cdot}17$ & $(0, -1)$ & $\text{yes}$ & $-51$ & $(48(5 - 1\sqrt{17}))^3(-33 +
8\sqrt{17})$\\
new & $(0, 1)$ & $\text{yes}$ & $-43$ & $-2^{18}~3^{3}~5^{3}$\\
 & $(\frac{1}{2}, \frac{-9}{8})$ & $\text{yes}$ & $-4$ & $2^{6}~3^{3}$\\
 & $(\frac{1}{2}, \frac{9}{8})$ & $\text{no}$ & $$ & $\Q(\sqrt{1081}, 
\sqrt{-779263})$\\
\midrule
230 & $\infty'$ & $\text{yes}$ & $-40$ & $(6(65 + 27\sqrt{5}))^3$\\
$2{\cdot}5{\cdot}23$ & $(0, -1)$ & $\text{yes}$ & $-15$ & $1(3(5 - 
4\sqrt{5}))^3(-3 + 1\sqrt{5})/2$\\
 & $(0, 1)$ & $\text{yes}$ & $-20$ & $(2(25 + 13\sqrt{5}))^3$\\
 & $(1, -5)$ & $\text{yes}$ & $-120$ & $\Q(\sqrt{2}, \sqrt{5})$\\
 & $(1, 5)$ & $\text{yes}$ & $-520$ & $\Q(\sqrt{5}, \sqrt{13})$\\
 & $(-2, -5)$ & $\text{yes}$ & $-1380$ & $\Q(\sqrt{3}, \sqrt{5}, \sqrt{23})$\\ \hdashline
&$(-2,  5) $&$\text{yes, yes}$&$-15,-60$&
\begin{tabular}{c}
$1(3(-5 + 4\sqrt{5}))^3(3 - 1\sqrt{5})/2$, \\
$1(3(470 + 213\sqrt{5}))^3(1 + 1\sqrt{5})/2$ \\
\end{tabular}\\ \hdashline
 & $(3, -35)$ & $\text{yes}$ & $-180$ & $\Q(\sqrt{3}, \sqrt{5})$\\
 & $(3, 35)$ & $\text{no}$ & $$ & $\Q(\sqrt{685}, \sqrt{19043}, \sqrt{19317})$\\
\midrule
255 & $\infty'$ & $\text{yes}$ & $-795$ & $\Q(\sqrt{5}, \sqrt{53})$\\
$3{\cdot}5{\cdot}17$ & $(-1, -2)$ & $\text{yes}$ & $-15$ & $1(3(5 - 
4\sqrt{5}))^3(-3 + 1\sqrt{5})/2$\\
new & $(-1, 2)$ & $\text{yes}$ & $-240$ & $\Q(\sqrt{3}, \sqrt{5})$\\
 & $(0, -1)$ & $\text{yes}$ & $-35$ & $(16(-15 + 7\sqrt{5}))^3$\\
 & $(0, 1)$ & $\text{yes}$ & $-51$ & $(48(-5 + 1\sqrt{17}))^3(33 - 
8\sqrt{17})$\\
 & $(-2, -5)$ & $\text{yes}$ & $-3315$ & $\Q(\sqrt{5}, \sqrt{13}, \sqrt{17})$\\
 & $(-2, 5)$ & $\text{yes}$ & $-120$ & $\Q(\sqrt{2}, \sqrt{5})$\\
 & $(\frac{1}{2}, \frac{-5}{8})$ & $\text{yes}$ & $-195$ & $\Q(\sqrt{5}, 
\sqrt{13})$\\
 & $(\frac{1}{2}, \frac{5}{8})$ & $\text{no}$ & $$ & $\Q(\sqrt{-15}, \sqrt{-87},
\sqrt{-3783})$\\
 & $(\frac{3}{5}, \frac{-38}{125})$ & $\text{no}$ & $$ & $\Q(\sqrt{19}, 
\sqrt{-6931}, \sqrt{1330349})$\\
 & $(\frac{3}{5}, \frac{38}{125})$ & $\text{yes}$ & $-60$ & $1(3(470 + 
213\sqrt{5}))^3(1 + 1\sqrt{5})/2$\\
\midrule
266 & $\infty'$ & $\text{yes}$ & $-52$ & $(30(31 + 9\sqrt{13}))^3$\\
$2{\cdot}7{\cdot}19$ & $(-1, -1)$ & $\text{yes, yes}$ & $-12, -3$ & 
$2^{4}~3^{3}~5^{3}, 0$\\
new & $(-1, 1)$ & $\text{yes}$ & $-84$ & $\Q(\sqrt{3}, \sqrt{7})$\\
 & $(0, -1)$ & $\text{yes}$ & $-40$ & $(6(65 + 27\sqrt{5}))^3$\\
 & $(0, 1)$ & $\text{yes}$ & $-280$ & $\Q(\sqrt{2}, \sqrt{5})$\\
 & $(\frac{-5}{2}, \frac{-83}{8})$ & $\text{yes}$ & $-532$ & $\Q(\sqrt{7}, 
\sqrt{19})$\\
 & $(\frac{-5}{2}, \frac{83}{8})$ & $\text{no}$ & $$ & $\Q(\sqrt{105}, 
\sqrt{-415}, \sqrt{1041})$\\
\midrule
285 & $\infty'$ & $\text{yes}$ & $-51$ & $(48(-5 + 1\sqrt{17}))^3(33 - 
8\sqrt{17})$\\
$3{\cdot}5{\cdot}19$ & $(-1, -4)$ & $\text{yes}$ & $-60$ & $1(3(470 + 
213\sqrt{5}))^3(1 + 1\sqrt{5})/2$\\
 & $(-1, 4)$ & $\text{yes}$ & $-15$ & $1(3(-5 + 4\sqrt{5}))^3(3 - 
1\sqrt{5})/2$\\
 & $(0, 0)$ & $\text{yes, yes}$ & $-3, -75$ & $0, (48(69 - 31\sqrt{5}))^3(0 + 
1\sqrt{5})$\\
 & $(3, -24)$ & $\text{yes}$ & $-240$ & $\Q(\sqrt{3}, \sqrt{5})$\\
 & $(3, 24)$ & $\text{no}$ & $$ & $\Q(\sqrt{3}, \sqrt{95}, \sqrt{60197})$\\
 & $(\frac{-3}{2}, \frac{-57}{8})$ & $\text{yes}$ & $-1995$ & $\Q(\sqrt{5}, 
\sqrt{21}, \sqrt{57})$\\
 & $(\frac{-3}{2}, \frac{57}{8})$ & $\text{no}$ & $$ & $\Q(\sqrt{57}, 
\sqrt{-79}, \sqrt{11985})$\\
\midrule
286 & $\infty'$ & $\text{yes}$ & $-40$ & $(6(65 + 27\sqrt{5}))^3$\\
$2{\cdot}11{\cdot}13$ & $(-1, -4)$ & $\text{yes}$ & $-88$ & $(60(155 + 
108\sqrt{2}))^3$\\
 & $(-1, 4)$ & $\text{yes}$ & $-52$ & $(30(31 + 9\sqrt{13}))^3$\\
 & $(\frac{5}{2}, \frac{-143}{8})$ & $\text{no}$ & $$ & $\Q(\sqrt{-407}, 
\sqrt{1841}, \sqrt{-3367})$\\
 & $(\frac{5}{2}, \frac{143}{8})$ & $\text{no}$ & $$ & $\Q(\sqrt{39}, 
\sqrt{231}, \sqrt{168917})$\\
\midrule
287 & $\infty'$ & $\text{yes}$ & $-91$ & $(48(-227 + 63\sqrt{13}))^3$\\
$7{\cdot}41$ & $(-2, -9)$ & $\text{no}$ & $$ & $\Q(\sqrt{8321}, 
\sqrt{349013})$\\
new & $(-2, 9)$ & $\text{yes}$ & $-1435$ & $\Q(\sqrt{5}, \sqrt{41})$\\
\midrule
299 & $\infty'$ & $\text{yes}$ & $-91$ & $(48(-227 + 63\sqrt{13}))^3$\\
$13{\cdot}23$ & $(\frac{-1}{2}, \frac{-1}{8})$ & $\text{no}$ & $$ & 
$\Q(\sqrt{1513}, \sqrt{-3325543})$\\
new & $(\frac{-1}{2}, \frac{1}{8})$ & $\text{yes}$ & $-43$ & 
$-2^{18}~3^{3}~5^{3}$\\
\midrule
    330 & $\infty'$ & yes & $-120$ & $\Q(\sqrt{2},\sqrt{5})$\\
    $2{\cdot}3{\cdot}5{\cdot}11$ & $(-1,-4)$ & yes & $-660$ & $\Q(\sqrt{3},\sqrt{5},\sqrt{11})$\\
    new & $(-1,4)$ & yes & $-24$ & $(12(5 + 2\sqrt{2}))^3(3 + 2\sqrt{2})$\\
     & $(0,-3)$ & yes & $-84$ & $\Q(\sqrt{3},\sqrt{7})$\\
     & $(0,3)$ & yes & $-420$ & $\Q(\sqrt{3},\sqrt{5},\sqrt{7})$\\
     & $(\frac{7}{8},\frac{-3497}{512})$ & no &  & $\Q(\sqrt{-95},\sqrt{3497},\sqrt{-20895},\sqrt{-33807})$\\
     & $(\frac{7}{8},\frac{3497}{512})$ & yes & $-1320$ & $\Q(\sqrt{5},\sqrt{6},\sqrt{22})$\\
\midrule
357 & $\infty'$ & $\text{yes}$ & $-168$ & $\Q(\sqrt{6}, \sqrt{14})$\\
$3{\cdot}7{\cdot}17$ & $(2, -14)$ & $\text{yes}$ & $-35$ & $(16(-15 + 
7\sqrt{5}))^3$\\
 & $(2, 14)$ & $\text{no}$ & $$ & $\Q(\sqrt{21}, \sqrt{293}, \sqrt{89997})$\\
\midrule
390 & $\infty'$ & $\text{yes}$ & $-5460$ & $\Q(\sqrt{3}, \sqrt{5}, \sqrt{7}, 
\sqrt{13})$\\
$2{\cdot}3{\cdot}5{\cdot}13$ & $(0, -1)$ & $\text{yes}$ & $-420$ & $\Q(\sqrt{3},
\sqrt{5}, \sqrt{7})$\\
 & $(0, 1)$ & $\text{yes}$ & $-120$ & $\Q(\sqrt{2}, \sqrt{5})$\\
 & $(1, -2)$ & $\text{yes, yes}$ & $-4, -36$ & $2^{6}~3^{3}, (4(21 + 
20\sqrt{3}))^3(7 + 4\sqrt{3})$\\
 & $(1, 2)$ & $\text{yes}$ & $-660$ & $\Q(\sqrt{3}, \sqrt{5}, \sqrt{11})$\\
\end{longtable}
}
\end{appendix}

\section*{Data Availability Statement}
The code to reproduce our computations, and the log files that were produced with Magma~2.26-10, are available in a public GitHub repository~\cite{ACKP}. Magma code to compute the data in Table~\ref{tab:jinv} and generate the \LaTeX{} source for the table, is also available in the same repository.

\bibliographystyle{amsalpha}
\bibliography{references}

\end{document}